\documentclass[
11pt,
draftclsnofoot,
onecolumn,
usletter,
]{IEEEtran}

\usepackage{amsfonts,amssymb,amsmath,amsthm,dsfont}

\usepackage{subfigure}
\usepackage{url, fullpage, color, upgreek}
\usepackage{graphicx}

\newcommand{\Rbb}{\mathbb{R}}
\newcommand{\scp}[2]{\langle #1, #2 \rangle}
\newcommand{\Nbb}{\mathbb{N}}
\newcommand{\Zbb}{\mathbb{Z}}
\newtheorem{theorem}{Theorem}
\newtheorem{definition}{Definition}
\newtheorem{proposition}{Proposition}
\newtheorem{lemma}{Lemma}

\newcommand{\inv}[1]{\frac{1}{#1}}
\newcommand{\supp}{{\rm supp}\,}
\newcommand{\tinv}[1]{{\textstyle\frac{1}{#1}}}
\newcommand{\sign}{{\rm sign}\,}
\newcommand{\ud}{\mathrm{d}} 
\renewcommand{\leq}{\leqslant}
\renewcommand{\geq}{\geqslant}
\newcommand{\E}{{\mathbb{E}}}
\newcommand{\Prob}{{\mathbb{P}}}

\newcommand{\ie}{\mbox{{\em i.e.},~}}
\newcommand{\eg}{\mbox{{\em e.g.},~}}

\DeclareMathOperator{\prox}{prox}
\newcommand{\Proj}{\mathcal{P}}
\DeclareMathOperator{\Id}{\mathds{1}}
\DeclareMathOperator{\dom}{dom}

\DeclareMathOperator*{\argmin}{argmin}
\newcommand{\norm}[1]{\|#1\|}
\newcommand{\Hm}{\mathcal{H}}

\title{Dequantizing Compressed Sensing:\\[2mm] \Large 
When Oversampling and Non-Gaussian Constraints Combine}

\author{L. Jacques, D. K. Hammond, M. J. Fadili\\[2mm]
\thanks{LJ is with the Information and Communication Technologies, Electronics
and Applied Mathematics (ICTEAM) Sector, Universit\'e catholique de Louvain
(UCL), Belgium. LJ is a Postdoctoral Researcher of the Belgian National Science Foundation
(F.R.S.-FNRS).}
\thanks{DKH is with the Neuroinformatics Center, University of Oregon, USA.}
\thanks{MJF is with the GREYC CNRS-ENSICAEN, Universit\'e de Caen, France.}
\thanks{A part of this work was presented at the IEEE Intl. Conf.
    Image Proc. (ICIP), Cairo, Egypt, 2009 \cite{jahafa09}.}  \vspace{-.5cm} }

\begin{document}

\maketitle

\begin{abstract}
  In this paper we study the problem of recovering sparse or
  compressible signals from uniformly quantized measurements.  We
  present a new class of convex optimization programs, or decoders,
  coined Basis Pursuit DeQuantizer of moment $p$ (BPDQ$_p$), that
  model the quantization distortion more faithfully than the commonly
  used Basis Pursuit DeNoise (BPDN) program. Our decoders proceed by
  minimizing the sparsity of the signal to be reconstructed subject to
  a data-fidelity constraint expressed in the $\ell_p$-norm
  of the residual error for $2\leq p\leq \infty$.

  We show theoretically that, \emph{(i)} the reconstruction error of
  these new decoders is bounded if the sensing matrix satisfies an
  extended Restricted Isometry Property involving the $\ell_p$ norm,
  and \emph{(ii)}, for Gaussian random matrices and uniformly
  quantized measurements, BPDQ$_p$ performance exceeds that of BPDN by
  dividing the reconstruction error due to quantization by
  $\sqrt{p+1}$. This last effect happens with high probability when
  the number of measurements exceeds a value growing with $p$, \ie in
  an oversampled situation compared to what is commonly required by
  BPDN = BPDQ$_2$. To demonstrate the theoretical power of BPDQ$_p$,
  we report numerical simulations on signal and image reconstruction
  problems.
\end{abstract}


\section{Introduction}
\label{sec:intro}

The theory of Compressed Sensing (CS)
\cite{candes2006qru,donoho2006cs} aims at reconstructing sparse or
compressible signals from a small number of linear measurements
compared to the dimensionality of the signal space.  In short, the
signal reconstruction is possible if the underlying sensing matrix is
well behaved, \ie if it respects a Restricted Isometry Property (RIP)
saying roughly that any small subset of its columns is ``close'' to an
orthogonal basis. The signal recovery is then obtained using
non-linear techniques based on convex optimization promoting signal
sparsity, such as the Basis Pursuit program \cite{donoho2006cs}. What
makes CS more than merely an interesting theoretical concept is that
some classes of randomly generated matrices (\eg Gaussian, Bernoulli,
partial Fourier ensemble, etc) satisfy the RIP with overwhelming
probability. This happens as soon as their number of rows, \ie the
number of CS measurements, is higher than a few multiples of the
assumed signal sparsity.

In a realistic acquisition system, quantization of these measurements
is a natural process that Compressed Sensing theory has to handle
conveniently. 
One commonly used technique is to simply treat the quantization
distortion as Gaussian noise, which leads to reconstruction based on
solving the Basis Pursuit DeNoising (BPDN) program (either in its
constrained or augmented Lagrangian forms) \cite{Candes2004}.  While
this approach can give acceptable results, it is theoretically
unsatisfactory as the measurement error created by quantization is
highly non-Gaussian, being essentially uniform and bounded by the
quantization bin width.

An appealing requirement for the design of better reconstruction
methods is the Quantization Consistency (QC) constraint, \ie that the
requantized measurements of the reconstructed signal equal the
original quantized measurements.  This idea, in some form, has
appeared previously in the literature.  Near the beginning of the
development of CS theory, Cand\`es et al. mentioned that the
$\ell_2$-norm of BPDN should be replaced by the $\ell_\infty$-norm to
handle more naturally the quantization distortion of the measurements
\cite{Candes2004}. More recently, in \cite{CS1bit}, the extreme case
of 1-bit CS is studied, \ie when only the signs of the measurements
are sent to the decoder. Authors tackle the reconstruction problem by
adding a sign consistency constraint in a modified BPDN program
working on the sphere of unit-norm signals. In \cite{Dai2009}, an
adaptation of both BPDN and the Subspace Pursuit integrates an
explicit QC constraint. In \cite{Zymnis2009}, a model integrating
additional Gaussian noise on the measurements before their
quantization is analyzed and solved with a $\ell_1$-regularized maximum
likelihood program. However, in spite of interesting experimental
results, no theoretical guarantees are given about the approximation
error reached by these solutions. The QC constraint has also been used
previously for image and signal processing outside of the CS field.
Examples include oversampled Analog to Digital Converters (ADC)
\cite{thao1994dao}, and in image restoration problems
\cite{weiss2006sac, weiss2008car}.

In this paper, we propose a new class of convex optimization programs,
or decoders, coined the Basis Pursuit DeQuantizer of moment $p$
(BPDQ$_p$) that model the quantization distortion more
faithfully. These proceed by minimizing the sparsity of the
reconstructed signal (expressed in the $\ell_1$-norm) subject to a
particular data-fidelity constraint. This constraint imposes that the
difference between the original and the reproduced measurements have
bounded $\ell_p$-norm, for $2\leq p\leq \infty$. As $p$ approaches
infinity, this fidelity term reproduces the QC constraint as promoted
initially in \cite{Candes2004}. However, our idea is to study, given a
certain sparsity level and in function of the number of measurements
available, which moment $2\leq p\leq \infty$ provides the best
reconstruction result.

Our overall result, which surprisingly does not favor $p=\infty$, may
be expressed by the principle: \emph{Given a certain sparsity level,
  if the number of measurements is higher than a minimal value growing
  with $p$, \ie in oversampled situations, by using BPDQ$_p$ instead
  of BPDN~=~BPDQ$_2$ the reconstruction error due to quantization can
  be reduced by a factor of $\sqrt{p+1}$.}


At first glance, it could seem counterintuitive to oversample the
``compressive sensing'' of a signal. After all, many results in
Compressed Sensing seek to limit the number of measurements required
to encode a signal, while guaranteeing exact reconstruction with high
probability. However, as analyzed for instance in
\cite{goyal2008compressive}, this way of thinking avoids to
considering the actual amount of information needed to describe the
measurement vector. In the case of noiseless observations of a sparse
signal, Compressed Sensing guarantees perfect reconstruction only for
real-valued measurements, \ie for an infinite number of bits per
measurements.

From a rate-distortion perspective, the analysis shown in
\cite{cand2006encoding,bouf2007quantspars} demonstrates also that CS
is suboptimal compared to transform coding. Under that point of view,
the best CS encoding strategy is to use all the available bit-rate to
obtain as few CS measurements as possible and quantize them as finely
as possible.

However, in many practical situations the quantization bit-depth per
measurement is pre-determined by the hardware, \eg for real sensors
embedding CS and a fixed A/D conversion of the measurements. In that
case, the only way to improve the reconstruction quality is to gather
more measurements, \ie to oversample the signal\footnote{Generally, it
  is also less expensive in hardware to oversample a signal than to
  quantize measurements more finely.}. This does not degrade one of
the main interests of Compressed Sensing, \ie providing highly
informative linear signal measurements at a very low computation cost.

The paper is structured as follows. In Section
\ref{sec:cs-quantiz-fmwk}, we review the principles of Compressed
Sensing and previous approaches for accommodating the problem of
measurement quantization. Section \ref{sec:new-class-decoders}
introduces the BPDQ$_p$ decoders. Their stability, \ie the
$\ell_2-\ell_1$ instance optimality, is deduced using an extended
version of the Restricted Isometry Property involving the
$\ell_p$-norm. In Section \ref{sec:example-rip_p-2}, Standard Gaussian
Random matrices, \ie whose entries are independent and identically
distributed (iid) standard Gaussian, are shown to satisfy this
property with high probability for a sufficiently large number of
measurements. Section \ref{sec:BPDQ-approx-error-for-quantiz} explains
the key result of this paper; that the approximation error of BPDQ$_p$
scales inversely with $\sqrt{p+1}$.  Section \ref{sec:implementations}
describes the convex optimization framework adopted to solve the
BPDQ$_p$ programs. Finally, Section \ref{sec:experiment} provides
experimental validation of the theoretical power of BPDQ$_p$ on 1-D
signals and on an image reconstruction example.

\section{Compressed Sensing and Quantization of Measurements}
\label{sec:cs-quantiz-fmwk}

In Compressed Sensing (CS) theory \cite{candes2006qru, donoho2006cs},
the signal $x\in\Rbb^N$ to be acquired and subsequently reconstructed
is typically assumed to be sparse or \emph{compressible} in an
orthogonal\footnote{A generalization for redundant basis, or
  dictionary, exists \cite{rauhut2008csa,ying2009lta}.} basis
$\Psi\in\Rbb^{N\times N}$ (\eg wavelet basis, Fourier, etc.).  In
other words, the best $K$-term approximation $x_K$ of $x$ in $\Psi$
gives an exact (for the sparse case) or accurate (for the compressible
case) representation of $x$ even for small $K<N$. For simplicity, only
the canonical basis $\Psi =\rm Id$ will be considered here.

At the acquisition stage, $x$ is encoded by $m$ linear measurements
(with $K\leq m\leq N$) provided by a sensing matrix
$\Phi\in\Rbb^{m\times N}$, \ie all known information about $x$ is
contained in the $m$ measurements $\scp{\varphi_i}{x}=\sum_k
\varphi^*_{ik} x_k$, where $\{\varphi_i\}_{i=0}^{m-1}$ are the rows of
$\Phi$.

In this paper, we are interested in a particular non-ideal sensing
model. Indeed, as measurement of continuous signals by digital devices
always involves some form of quantization, in practice devices based
on CS encoding must be able to accommodate the distortions in the
linear measurements created by quantization. Therefore, we adopt the
noiseless and uniformly quantized sensing (or coding) model:
\begin{equation}
  \label{eq:quantiz-fmwk}
  y_{\rm q}\ =\ Q_\upalpha[\Phi x]\ =\ \Phi x + n,
\end{equation}
where $y_{\rm q}\in(\upalpha\Zbb+\frac{\upalpha}{2})^m$ is the
quantized measurement vector, $(Q_\upalpha[\cdot])_i=\upalpha\lfloor
(\cdot)_i/\upalpha\rfloor + \tfrac{\upalpha}{2}$ is the uniform
quantization operator in $\Rbb^m$ of bin width $\upalpha$, and
$n\triangleq Q_\upalpha[\Phi x]\- \Phi x$ is the \emph{quantization
  distortion}.

The model (\ref{eq:quantiz-fmwk}) is a realistic description of
systems where the quantization distortion dominates other secondary
noise sources (\eg thermal noise), an assumption valid for many
electronic measurement devices including ADC. In this paper we
restrict our study to using this extremely simple uniform quantization
model, in order to concentrate on the interaction with the CS
theory. For instance, this quantization scenario does not take into
account the possible \emph{saturation} of the quantizer happening when
the value to be digitized is outside the operating range of the
quantizer, this range being determined by the number of bits
available. For Compressed Sensing, this effect has been studied
recently in \cite{laska09}. Authors obtained better reconstruction
methods by either imposing to reproduce saturated measurements
(Saturation Consistency) or by discarding these thanks to the
``democratic'' property of most of the random sensing matrices. Their
work however does not integrate the Quantization Consistency for all
the unsaturated measurements. The study of more realistic non-uniform
quantization is also deferred as a question for future research.

In much previous work in CS, the reconstruction of $x$
from $y_{\rm q}$ is obtained by  treating the
quantization distortion $n$ as a noise of bounded power
(\ie $\ell_2$-norm) $\|n\|^2_2=\sum_k |n_k|^2$.  
In this case, a
robust reconstruction of the signal $x$ from corrupted measurements
$y=\Phi x + n$ is provided by the Basis Pursuit DeNoise (BPDN) program
(or decoder) \cite{candes2006ssr}:
\begin{equation*}
  \Delta(y,\epsilon) = \argmin_{u\in\Rbb^N} \|u\|_1
  \ {\rm
  s.t.}\ \|y - \Phi u\|_2\leq \epsilon. \eqno{\textrm{\small ({\bf  BPDN})}}
\end{equation*}
This convex optimization program can be solved numerically by methods
like Second Order Cone Programming or by monotone operator splitting methods
\cite{Fadili2009,combettes2004smi} described in Section
\ref{sec:implementations}. Notice that the noiseless situation
$\epsilon=0$ leads to the Basis Pursuit (BP) program, which may also
be solved by Linear Programming \cite{Chen2001}.

An important condition for BPDN to provide a good reconstruction is the
\emph{feasibility} of the initial signal $x$, \ie we must chose
$\epsilon$ in the (\emph{fidelity}) constraint of BPDN such that
$\|n\|_2=\|y - \Phi x\|_2\leq \epsilon$. In \cite{candes2006ssr}, an
estimator of $\epsilon$ for $y=y_{\rm q}$ is obtained by considering
$n$ as a random vector $\xi\in\Rbb^m$ distributed uniformly over the
quantization bins, \ie $\xi_i \sim_{\rm iid}
U([-\tfrac{\upalpha}{2},\tfrac{\upalpha}{2}])$. 

An easy computation shows then that $\|\xi\|^2_2 \leq
\epsilon^2_2(\upalpha)$ with probability higher than
$1-e^{-c_0\kappa^2}$ for a certain constant $c_0>0$ (by the
Chernoff-Hoeffding bound \cite{hoeffding1963pis}), where
\begin{equation*}
\epsilon^2_2(\upalpha) \triangleq \E\|\xi\|^2_2 +
\kappa\,\sqrt{{\rm Var}\|\xi\|^2_2} = \tfrac{\upalpha^2}{12}m +
\kappa\,\tfrac{\upalpha^2}{6\sqrt{5}}\,m^{\frac{1}{2}}.
\end{equation*}

Therefore, CS usually handles quantization
distortion by setting $\epsilon=\epsilon_2(\upalpha)$, typically for
$\kappa=2$.

When the feasibility is satisfied, the stability of BPDN is guaranteed
if the sensing matrix $\Phi\in\Rbb^{m\times N}$ satisfies one instance
of the following property:
\begin{definition}
  A matrix $\Phi\in\Rbb^{m\times N}$ satisfies the (extended)
  Restricted Isometry Property (RIP$_{p,q}$) (with $p,q>0$) of order
  $K$ and radius $\delta_K\in (0,1)$, if there exists a constant
  $\mu_{p,q}>0$ such that
\begin{equation}
\label{eq:rippq}
\mu_{p,q}\,(1-\delta_K)^{1/q}\,\|u\|_q\, \leq\, \|\Phi u\|_p\, \leq\, \mu_{p,q}\,(1+\delta_K)^{1/q}\,\|u\|_q,
\end{equation}
for all $K$-sparse signals $u\in\Rbb^N$. 
\end{definition}
In other words, $\Phi$, as a mapping from $\ell^m_p=(\Rbb^m,
\|\!\cdot\!\|_p)$ to $\ell^N_q=(\Rbb^N, \|\!\cdot\!\|_q)$, acts as a
(scaled) isometry on $K$-sparse signals of $\Rbb^N$. This definition is more
general than the common RIP \cite{candes2008rip}. This latter, which
ensures the stability of BPDN (see Theorem
\ref{prop:l2-l1-instance-optimality-BP} below), corresponds to $p=q=2$
in (\ref{eq:rippq}). The original definition considers also normalized
matrices $\bar \Phi = \Phi /\mu_{2,2}$ having unit-norm columns (in
expectation) so that $\mu_{2,2}$ is absorbed in the normalizing
constant.

We prefer to use this extended RIP$_{p,q}$ since, as it will become
clear in Section \ref{sec:BPDQ-approx-error-for-quantiz}, the case
$p\geq 2$ and ${q=2}$ provides us the interesting embedding
(\ref{eq:rippq}) for measurement vectors corrupted by generalized
Gaussian and uniform noises. As explained below, this definition
includes also other RIP
generalizations~\cite{Berinde2008,chartrand2008restricted}.

We note that there are several examples already described in the
literature of classes of matrices which satisfy the RIP$_{p,q}$ for
specific values of $p$ and $q$. For instance, for $p=q=2$, a matrix
$\Phi\in\Rbb^{m\times N}$ with each of its entries drawn independently
from a (sub) Gaussian random variable satisfies this property with an
overwhelming probability if $m \geq c K \log N/K$ for some value $c>0$
independent of the involved dimensions \cite{JLmeetCS, DoTa09,
  mendelson2007reconstruction}.  This is the case of Standard Gaussian
Random (SGR) matrices whose entries are iid $\Phi_{ij}\,\sim\,
\mathcal{N}(0,1)$, and of the Bernoulli matrices with $\Phi_{ij} = \pm
1$ with equal probability, both cases having $\mu_{2,2}=\sqrt{m}$
\cite{JLmeetCS}.  Other random constructions satisfying the
RIP$_{2,2}$ are known (\eg partial Fourier ensemble)
\cite{candes2006qru,candes2006ssr}. For the case $p=q=1+O(1)/\log N$,
it is proved in \cite{Berinde2008,BerInd07} that sparse matrices
obtained from an adjacency matrix of a high-quality unbalanced
expander graph are RIP$_{p,p}$ (with $\mu_{p,p}^2=1/(1-\delta_K)$). In
the context of non-convex signal reconstruction, the authors in
\cite{chartrand2008restricted} show also that Gaussian random matrices
satisfy the Restricted $p$-Isometry, \ie RIP$_{p,q}$ for $q=2$,
$0<p<1$, $\mu_{p,2}=1$ and appropriate redefinition of $\delta_K$.
\medskip

The following theorem expresses the announced stability result,
\ie the $\ell_2-\ell_1$ instance optimality\footnote{Adopting the
  definition of mixed-norm instance optimality
  \cite{Cohen-bestkterm}.} of BPDN, as a consequence of the
RIP$_{2,2}$.
\begin{theorem}[\cite{candes2008rip}]
  \label{prop:l2-l1-instance-optimality-BP}
  Let $x\in\Rbb^N$ be a signal whose compressibility is measured by
  the decreasing of the $K$-term $\ell_1$-approximation error $e_0(K)
  = K^{-\frac{1}{2}}\,\|x - x_K\|_1$, for $0\leq K\leq N$, and $x_K$
  the best $K$-term $\ell_2$-approximation of $x$. Let $\Phi$ be a
  RIP$_{2,2}$ matrix of order $2K$ and radius
  $0<\delta_{2K}<\sqrt{2}-1$. Given a measurement vector $y=\Phi x +
  n$ corrupted by a noise $n$ with power $\|n\|_2\leq \epsilon$, the
  solution $x^*=\Delta(y,\epsilon)$ obeys  
\begin{equation}
\label{eq:l2-l1-instance-optimality-BP}
\|x^* - x\|_2\ \leq\ A\,e_0(K)\ +\ B\,\tfrac{\epsilon}{\mu_{2,2}},
\end{equation}
for  $A(\Phi,K) = 2\,\tfrac{1 + (\sqrt{2}-1)\delta_{2K}}{1 -
  (\sqrt{2}+1)\delta_{2K}}$ and $B(\Phi,K) = \tfrac{4\sqrt{1 + \delta_{2K}}}{1
  - (\sqrt{2}+1)\delta_{2K}}$. For instance, for $\delta_{2K}=0.2$,
$A<4.2$ and $B<8.5$.
\end{theorem}

Let us precise that the theorem condition $\delta_{2K}<\sqrt{2}-1$ on the RIP
radius can be refined (like in \cite{foucart2009sparsest}). We know
nevertheless from Davies and Gribonval \cite{davies2009ripfail} that
$\ell_1$-minimization will fail for at least one vector for
$\delta_{2K} > 1/\sqrt 2$. The room for improvement is then very
small.

\medskip

Using the BPDN decoder to account for quantization distortion is
theoretically unsatisfying for several reasons. First, there is no
guarantee that the BPDN solution $x^*$ respects the Quantization
Consistency, \ie
\begin{equation*}
  Q_\upalpha[\Phi x^*]=y_{\rm q}\ \Leftrightarrow\ \|y_{\rm q} - \Phi
  x^*\|_\infty \leq \tfrac{\upalpha}{2},\eqno{({\bf QC})}
\end{equation*}
which is not necessarily implied by the BPDN $\ell_2$ fidelity
constraint. The failure of BPDN to respect QC suggests that it may not
be taking advantage of all of the available information about the
noise structure in the measurements.

Second, from a Bayesian Maximum a Posteriori (MAP) standpoint, BPDN
can be viewed as solving an ill-posed inverse problem where the
$\ell_2$-norm used in the fidelity term corresponds to the conditional
log-likelihood associated to an additive white Gaussian
noise. However, the quantization distortion is not Gaussian, but
rather uniformly distributed. This motivates the need for a new kind
of CS decoder that more faithfully models the quantization distortion.

\section{Basis Pursuit DeQuantizer (BPDQ$_p$)}
\label{sec:def-bpdq}
\label{sec:new-class-decoders}


The considerations of the previous section encourage the definition of
a new class of optimization programs (or decoders) generalizing the
fidelity term of the BPDN program.

Our approach is based on reconstructing a sparse approximation of $x$
from its measurements $y = \Phi x + n$ under the assumption that
$\ell_p$-norm ($p\geq 1$) of the noise $n$ is bounded, \ie
$\|n\|_p^p=\sum_k|n_k|^p\leq \epsilon^p$ for some $\epsilon>0$.  We
introduce the novel programs
\begin{equation*} \label{eq:BPDQ}
  \Delta_{p}(y,\epsilon) = \argmin_{u\in\Rbb^N}
  \|u\|_1\ {\rm s.t.}\ \|y - \Phi u\|_p \leq \epsilon. \eqno{\textrm{\small ({\bf  BPDQ}$_p$)}}
\end{equation*}
The fidelity constraint expressed in the $\ell_p$-norm is now tuned to
noises that follow a zero-mean Generalized Gaussian
Distribution\footnote{The probability density function $f$ of such a
  distribution is $f(x)\propto \exp (- |x/b|^p)$ for a standard
  deviation $\sigma \propto b$.} (GGD) of \emph{shape parameter} $p$
\cite{varanasi1989pgg}, with the uniform noise case corresponding to
$p\to\infty$.

We dub this class of decoders \emph{Basis Pursuit DeQuantizer} of
\emph{moment} $p$ (or BPDQ$_p$) since, for reasons that will become
clear in Section \ref{sec:BPDQ-approx-error-for-quantiz}, their
approximation error when $\Phi x$ is uniformly quantized has an
interesting decreasing behavior when both the moment $p$ and the
oversampling factor $m/K$ increase. Notice that the decoder
corresponding to $p=1$ has been previously analyzed in
\cite{fuchs2009fast} for Laplacian noise.


\medskip

One of the main results of this paper concerns the $\ell_2-\ell_1$
instance optimality of the BPDQ$_p$ decoders, \ie their stability
when the signal to be recovered is compressible, and when the
measurements are contaminated by noise of bounded $\ell_p$-norm. In
the following theorem, we show that such an optimality happens when
the sensing matrix respects the (extended) Restricted Isometry
Property RIP$_{p,2}$ for $2\leq p<\infty$.

\begin{theorem}
  \label{prop:l2-l1-instance-optimality-BPDQp}
  Let $x\in\Rbb^N$ be a signal with a
  $K$-term $\ell_1$-approximation error $e_0(K) =
  K^{-\frac{1}{2}}\,\|x - x_K\|_1$, for $0\leq K\leq N$ and $x_K$ the
  best $K$-term $\ell_2$-approximation of $x$. Let $\Phi$ be a
  RIP$_{p,2}$ matrix on $s$ sparse signals with constants
  $\delta_{s}$, for $s\in\{K,2K,3K\}$ and $2\leq p <\infty$. Given a
  measurement vector $y=\Phi x + n$ corrupted by a noise $n$ with
  bounded $\ell_p$-norm, \ie $\|n\|_p\leq \epsilon$, the solution
  $x^*_p=\Delta_{p}(y,\epsilon)$ of BPDQ$_p$ obeys
\begin{equation*}
\|x_p^* - x\|_2\ \leq\ A_{p}\,e_0(K)\ +\ B_{p}\,\epsilon/\mu_{p,2},
\end{equation*}
for values $A_p(\Phi,K) = \tfrac{2(1+C_p-\delta_{2K})}{1-\delta_{2K}-C_p}$,
$B_p(\Phi,K) = \tfrac{4\sqrt{1+\delta_{2K}}}{1-\delta_{2K}-C_p}$, and
$C_p=C_p(\Phi,2K,K)$ given in the proof of Lemma \ref{lemma:bound-scp-lp}
(Appendix \ref{sec:proof-lemma-bound-scp-lp}).
\end{theorem}

As shown in Appendix \ref{sec:proof-theor-refpr}, this theorem follows
from a generalization of the fundamental result proved by Cand\`es
\cite{candes2008rip} to the particular geometry of Banach spaces~$\ell_p$.

\section{Example of RIP$_{p,2}$ Matrices}
\label{sec:example-rip_p-2}

Interestingly, it turns out that SGR matrices $\Phi\in\Rbb^{m\times
  N}$ also satisfy the RIP$_{p,2}$ with high probability provided that
$m$ is sufficiently large compared to the sparsity $K$ of the signals
to measure. This is made formal in the following Proposition, for
which the proof\footnote{Interestingly, this proof shows also that SGR
  matrices are RIP$_{p,2}$ with high probability for $1<p< 2$ when $m$
  exceeds a similar bound to (\ref{eq:SGR-RIP-measur-bound}).} is
given in Appendix~\ref{app:why-rip-p}.

\begin{proposition}
\label{prop:gauss-rip-inf}
Let $\Phi\in\Rbb^{m\times N}$ be a Standard Gaussian Random (SGR)
matrix, \ie its entries are iid $\mathcal{N}(0,1)$. Then, if $m\geq
(p-1)2^{p+1}$ for $2\leq p<\infty$ and $m\geq 0$ for $p=\infty$, there
exists a constant $c>0$ such that, for
\begin{equation}
\label{eq:SGR-RIP-measur-bound}
\Theta_p(m)\,\geq\,c\,\delta^{-2}\,\big(K \log[e\tfrac{N}{K}(1
+ 12\delta^{-1})] + \log \tfrac{2}{\eta}\big),
\end{equation}
with $\Theta_p(m)=m^{2/p}$ for $1\leq p<\infty$ and $\Theta_p(m)=\log m$
for $p=\infty$, $\Phi$ is RIP$_{p,2}$ of order $K$ and radius $\delta$
with probability higher than $1-\eta$. Moreover, the value
$\mu_{p,2}=\E\|\xi\|_p$ is the expectation value of the
$\ell_p$-norm of a SGR vector $\xi\in\Rbb^m$.
\end{proposition}

Roughly speaking, this proposition tells us that to generate a matrix
that is RIP$_{p,2}$ with high probability, we need a number of
measurements $m$ that grows polynomially in $K\log N/K$ with an
``order'' $p/2$ for $2\leq p<\infty$, while the limit case $p=\infty$
grows exponentially in $K\log N/K$.

Notice that an asymptotic estimation of $\mu_{p,2}$, \ie for $m
\to\infty$, can be found in \cite{franccois2007cfd} for $1\leq
p<\infty$. However, as presented in the following Lemma (proved in
Appendix \ref{sec:proof-lemma-refl}), non-asymptotic bounds for
$\mu_{p,2}=\E\|\xi\|_p$ can be expressed in terms of
$$
(\E\|\xi\|^p_p)^{1/p}= (m\,\E|g|^p)^{1/p} = \nu_p\, m^{1/p}, 
$$
with $g\sim N(0,1)$ and $\nu_p^p = \E|g|^p =
2^{\frac{p}{2}}\,\pi^{-\inv{2}}\,\Gamma(\tfrac{p+1}{2})$.
\begin{lemma}
\label{lem:strict-bounds-mu_p}
If $\xi\in\Rbb^m$ is a SGR vector, then, for $1\leq p<\infty$, 
\begin{equation*}
  \big(1\ +\
  \tfrac{2^{p+1}}{m}\big)^{\frac{1}{p}-1}\,(\E\|\xi\|_p^p)^{\inv{p}}\ 
  \leq\ \E\|\xi\|_p\ \leq\ (\E\|\xi\|_p^p)^{\inv{p}}.
\end{equation*}
In particular, as soon as $m\geq \beta^{-1}\,2^{p+1}$ for $\beta\geq
0$, $ \E\|\xi\|_p \geq (\E\|\xi\|_p^p)^{\inv{p}}\,(1 +
\beta)^{\frac{1}{p}-1} \geq (\E\|\xi\|_p^p)^{\inv{p}}\,(1 -
\tfrac{p-1}{p}\beta)$. For $p=\infty$, there exists a
$\rho>0$ such that $\rho^{-1}\,\sqrt{\log m} \leq \E\|\xi\|_\infty
\leq \rho\,\sqrt{\log m}$.
\end{lemma}

An interesting aspect of matrices respecting the RIP$_{p,2}$ is that they
approximately preserve the decorrelation of sparse vectors
of disjoint supports. 
\begin{lemma}
\label{lemma:bound-scp-lp}
Let $u,v\in\Rbb^N$ with $\|u\|_0=s$ and $\|v\|_0=s'$ and\ \ $\supp\!(u)\ \cap\
\supp\!(v) = \emptyset$, and $2\leq p < \infty$. If $\Phi$ is RIP$_{p,2}$ of
order $s+s'$ with constant $\delta_{s+s'}$, and of orders $s$ and $s'$
with constants $\delta_s$ and $\delta_{s'}$, then
\begin{equation}
\label{eq:decor-J-phi-preserv}
|\scp{J(\Phi u)}{\Phi v}|\ \leq\ \mu_{p,2}^2\,C_p\,\|u\|_2\|v\|_2,
\end{equation}
with $(J(u))_i = \|u\|_p^{2-p}\,|u_i|^{p-1}\,\sign u_i$ and $C_p=C_p(\Phi,s,s')$ is
given explicitly in Appendix \ref{sec:proof-lemma-bound-scp-lp}.
\end{lemma}

It is worth mentioning that the value $C_p$ behaves as
$\sqrt{(\delta_{s}+\delta_{s+s'})\,(1+\delta_{s'})\,(p-2)}$ for large
$p$, and as $\delta_{s+s'} + \tfrac{3}{4}(1+\delta_{s+s'})(p-2)$ for
$p\simeq 2$.  Therefore, this result may be seen as a generalization
of the one proved in \cite{candes2008rip} (see Lemma 2.1) for $p=2$
with $C_2=\delta_{s+s'}$. As shown in Appendix
\ref{sec:proof-lemma-bound-scp-lp}, this Lemma uses explicitly the
2-smoothness of the Banach spaces $\ell_p$ when $p\geq 2$
\cite{bynum1976wpl,xu1991ibs}, in connection with the \emph{normalized
  duality mapping} $J$ that plays a central role in the geometrical
description of $\ell_p$.

Lemma \ref{lemma:bound-scp-lp} is at the heart of the proof of Theorem
\ref{prop:l2-l1-instance-optimality-BPDQp}, which prevents the later
from being valid for $p=\infty$. This is related to the fact that the
$\ell_{\infty}$ Banach space is not 2-smooth and no duality mapping
exists.  Therefore, any result for $p=\infty$ would require different
tools than those developed here.

\section{BPDQ$_p$ and Quantization Error Reduction}
\label{sec:BPDQ-approx-error-for-quantiz}

Let us now observe the particular behavior of the BPDQ$_p$ decoders on
quantized measurements of a sparse or compressible signal assuming
that $\upalpha$ is known at the decoding step. In this Section, we
consider that $p\geq 2$ everywhere.

First, if we assume in the model (\ref{eq:quantiz-fmwk}) that the
quantization distortion $n = Q_\upalpha[\Phi x] - \Phi x$ is uniformly
distributed in each quantization bin, the simple Lemma below provides
precise estimator $\epsilon$ for any $\ell_p$-norm of~$n$.

\begin{lemma}
\label{lemma:expec-and-val-lp-norm-unif-vec}
If $\xi\in\Rbb^m$ is a uniform random vector with $\xi_i\sim_{\rm iid}
U([-\tfrac{\upalpha}{2},\tfrac{\upalpha}{2}])$, then, for $1\leq p < \infty$,
\begin{equation}
\label{eq:expec-and-val-lp-norm-unif-vec}
\zeta_p = \E\|\xi\|^p_p\ =\ \tfrac{\upalpha^p}{2^p(p+1)}\,m.
\end{equation}
In addition, for any $\kappa>0$, $\Prob\big[\|\xi\|^p_p \geq \zeta_p +
\kappa\,\tfrac{\upalpha^p}{2^p}\,\sqrt{m}\,\big] \ \leq\
e^{-2\kappa^2}$, while, $\lim_{p\to\infty} (\zeta_p + 
\kappa\tfrac{\upalpha^p}{2^p}\sqrt{m})^{\frac{1}{p}} = \tfrac{\upalpha}{2}$.
\end{lemma}
\noindent The proof is given in Appendix
\ref{sec:proof-expec-and-val-lp-norm-unif-vec}.

According to this result, we may set the $\ell_p$-norm
bound $\epsilon$ of the program BPDQ$_p$ to
\begin{equation}
\label{eq:epsilon-p-def}
\epsilon\ =\ \epsilon_p(\upalpha)\ \triangleq\  
\tfrac{\upalpha}{2\,(p+1)^{1/p}}\,
\big(\,m + \kappa\,(p+1)\,\sqrt{m}\,\big)^{\inv{p}},
\end{equation}
so that, for $\kappa = 2$, we know that $x$ is a feasible solution of the
BPDQ$_p$ fidelity constraint with a probability exceeding $1-e^{-8} > 1-
3.4 \times 10^{-4}$.  

Second, Theorem~\ref{prop:l2-l1-instance-optimality-BPDQp} points out
that, when $\Phi$ is RIP$_{p,2}$ with $2\leq p < \infty$, the
approximation error of the BPDQ$_p$ decoders is the sum of two terms:
one that expresses the \emph{compressibility error} as measured by
$e_0(K)$, and one, the \emph{noise error}, proportional to the ratio
$\epsilon/\mu_{p,2}$. In particular, by Lemma
\ref{prop:gauss-rip-inf}, for $m$ respecting
(\ref{eq:SGR-RIP-measur-bound}), a SGR sensing matrix of $m$ rows
induces with a controlled probability
\begin{equation}
\label{eq:l2l1-for-quantiz}
\|x - x^*_p\|_2\ \leq\ A_p\,e_0(K)\ +\ B_p
\frac{\epsilon_p(\upalpha)}{\mu_{p,2}}.
\end{equation}

Combining (\ref{eq:epsilon-p-def}) and the result of Lemma
\ref{lem:strict-bounds-mu_p}, we may bound the noise error for uniform
quantization more precisely. Indeed, for $2\leq p <\infty$, if $m\geq
(p-1)2^{p+1}$, $\mu_{p,2} \geq \tfrac{p-1}{p}\, \nu_p m^{\frac{1}{p}}$
with $\nu_p =
\sqrt{2}\,\pi^{-\frac{1}{2p}}\,\Gamma(\tfrac{p+1}{2})^{\frac{1}{p}}$.

In addition, using a variant of the Stirling formula found in
\cite{spira1971cgf}, we know that $|\Gamma(x) -
(\tfrac{2\pi}{x})^{\frac{1}{2}}(\tfrac{x}{e})^x|\leq \tinv{9
  x}\,(\tfrac{2\pi}{x})^{\frac{1}{2}}(\tfrac{x}{e})^x$ for $x\geq 1$.
Therefore, we compute easily that, for $x=(p+1)/2> 1$, $\nu_p \geq c^{1/p}\,
(\tfrac{p+1}{e})^{1/2}\geq c\, (\tfrac{p+1}{e})^{1/2}$ with $c =
\frac{8\sqrt 2}{9 \sqrt e} < 1$. Finally, by (\ref{eq:epsilon-p-def}), we see that,
\begin{multline}
\label{eq:frac-epsilon-mu-bound}
\frac{\epsilon_p(\upalpha)}{\mu_{p,2}} \leq
\tfrac{p}{p-1}\,\tfrac{9e}{16\sqrt 2}\,\big(\tinv{p+1}+\kappa\,\tfrac{1}{\sqrt{m}}\big)^{1/p}\ \frac{\upalpha}{\sqrt{p+1}}\\
<\ C \frac{\upalpha}{\sqrt{p+1}}, 
\end{multline}
with $C={9e}{/(8\sqrt 2)}<2.17$, where we used the bound
$\frac{p}{p-1}\leq 2$ and the fact that
$(\tinv{p+1}+\kappa\,\tfrac{1}{\sqrt{m}})^{1/p}<1$ if $m >
(\frac{p+1}{p}\kappa)^2=O(1)$.

In summary, we can formulate the following principle.

\noindent \textbf{Oversampling Principle.} \emph{The noise error term
  in the $\ell_2-\ell_1$ instance optimality relation
  (\ref{eq:l2l1-for-quantiz}) in the case of uniform quantization of
  the measurements of a sparse or compressible signal is divided by
  $\sqrt{p+1}$ in oversampled SGR sensing, \ie when the
  \emph{oversampling factor} $m/K$ is higher than a minimal value
  increasing with $p$.}

Interestingly, this follows the improvement achieved by adding a QC
constraint in the decoding of oversampled ADC signal conversion
\cite{thao1994dao}.

The oversampling principle requires some additional
explanations. Taking a SGR matrix, by Proposition
\ref{prop:gauss-rip-inf}, if $m_p$ is the smallest number of
measurements for which such a randomly generated matrix $\Phi$ is
RIP$_{p,2}$ of radius $\delta_p<1$ with a certain nonzero probability,
taking $m>m_p$ allows one to generate a new random matrix with a
smaller radius $\delta < \delta_p$ with the same probability of
success.

Therefore, increasing the \emph{oversampling factor} $m/K$ provides
two effects. First, it enables one to hope for a matrix $\Phi$ that is
RIP$_{p,2}$ for high $p$, providing the desired error division by
$\sqrt{p+1}$. Second, as shown in Appendix \ref{sec:delta-propto-m1p},
since $\delta = O(m^{-1/p}\sqrt{\log m})$, oversampling gives a
smaller $\delta$ hence counteracting the increase of $p$ in the factor
$C_p$ of the values $A_p\geq 2$ and $B_p\geq 4$. This decrease of
$\delta$ also favors BPDN, but since the values $A=A_2\geq 2$ and
$B=B_2\geq 4$ in (\ref{eq:l2-l1-instance-optimality-BP}) are also
bounded from below this effect is limited. Consequently, as the
number of measurements increases the improvement in reconstruction
error for BPDN will saturate, while for BPDQ$_p$ the error will be
divided by $\sqrt{p+1}$.

From this result, it is very tempting to choose an extremely large
value for $p$ in order to decrease the noise error term
\eqref{eq:l2l1-for-quantiz}. There are however two obstacles with
this. First, the instance optimality result of
Theorem~\ref{prop:l2-l1-instance-optimality-BPDQp} is not directly
valid for $p=\infty$. Second, and more significantly, the necessity of
satisfying RIP$_{p,2}$ implies that we cannot take $p$ arbitrarily large
in Proposition~\ref{prop:gauss-rip-inf}. Indeed, for a given
oversampling factor $m/K$, a SGR matrix $\Phi$ can be RIP$_{p,2}$ only
over a finite interval $p\in [2, p_{\max}]$. This implies that for
each particular reconstruction problem, there should be an optimal
maximum value for $p$. We will demonstrate this effect experimentally
in Section \ref{sec:experiment}.

We remark that the compressibility error is not significantly reduced
by increasing $p$ when the number of measurements is large. This makes
sense as the $\ell_p$-norm appears only in the fidelity term of the
decoders, and we know that in the case where $\epsilon=0$ the
compressibility error remains in the BP decoder
\cite{candes2008rip}. Finally, note that due to the embedding of the
$\ell_p$-norms, \ie $\|\!\cdot\!\|_p\leq \|\!\cdot\!\|_{p'}$ if $p\geq
p'\geq 1$, increasing $p$ until $p_{\max}$ makes the fidelity term
closer to the QC.

\section{Numerical Implementation}
\label{sec:implementations}

This section is devoted to the description of the convex optimization
tools needed to numerically solve the Basis Pursuit DeQuantizer
program. While we generally utilize $p\geq 2$, the BPDQ$_p$ program is convex
for $p\geq 1$. In fact, the efficient iterative procedure we describe
will converge to to the global minimum of the BPDQ$_p$ program for all
$p\geq 1$.

\subsection{Proximal Optimization}
\label{sec:prox-optim}

The BPDQ$_p$ (and BPDN) decoders are special case of a general class of
convex problems \cite{Fadili2009, combettes2008pdm}
\begin{equation*}
  \label{eq:convex-prob}
  \arg\min_{x \in \Hm}\ f_1(x) + f_2(x),\eqno{({\rm \bf P})}
\end{equation*}
where $\Hm=\Rbb^N$ is seen as an Hilbert space equipped with the inner
product $\scp{x}{z}=\sum_{i} x_i z_i$. We denote by $\dom f
=\{x\in\Hm: f(x) < +\infty\}$ the domain of any $f:\Hm\to \Rbb \cup \{+\infty\}$. In
({$\rm \bf P$}), the functions $f_1, f_2:\Hm\to \Rbb \cup \{+\infty\}$ are assumed
\emph{(i)} convex functions which are not infinite everywhere,
\ie $\dom f_1,\dom f_2\neq \emptyset$, \emph{(ii)} $\dom
f_1\,\cap\,\dom f_2 \neq \emptyset$, and \emph{(iii)} these functions
are lower semi-continuous (lsc) meaning that $\liminf_{x
  \to x_0} f(x) = f(x_0)$ for all $x_0 \in \dom f$. The class of
functions satisfying these three properties is denoted
$\Gamma_0(\Rbb^N)$. For BPDQ$_p$, these two non-differentiable
functions are $f_1(x)=\|x\|_1$ and 
$f_2(x) = \imath_{T^p(\epsilon)}(x)= 0$ 
if $x\in T^p(\epsilon)$ and $\infty$ otherwise, \ie the indicator
function of the closed convex set 
$T^p(\epsilon)=\{x\in\Rbb^N:\|y_{\rm q}- \Phi x\|_p\leq \epsilon\}$.

It can be shown that the solutions of problem ({$\rm \bf P$}) are
characterized by the following fixed point equation: x solves 
({$\rm \bf P$}) if and only if
\begin{equation} 
\label{eq-fixedpoint}
x\ =\ (\Id\ +\ \beta\partial(f_1+f_2))^{-1}(x),\mbox{ \quad for }\beta > 0.
\end{equation}
The operator $\mathcal{J}_{\beta\partial{f}} = (\Id + \beta \partial f)^{-1}$ is
called the {\textit{resolvent operator}} associated to the
\emph{subdifferential operator} $\partial f$, $\beta$ is a positive
scalar known as the proximal step size, and $\Id$ is the identity map
on $\Hm$. We recall that the subdifferential of a function $f \in
\Gamma_0(\Hm)$ at $x \in \Hm$ is the set-valued map $
\partial f(x)\ =\ \{u \in \Hm: \forall z \in \Hm, f(z) \geq f(x) +
\scp{u}{z-x}\}$, where each element $u$ of $\partial f$ is called a
subgradient.

The resolvent operator is actually identified with the \emph{proximity
  operator} of $\beta f$, \ie $\mathcal{J}_{\beta\partial{f}} = \prox_{\beta
  f}$, introduced in \cite{moreau1962fcd} as a generalization of
convex projection operator. It is defined as the unique solution
$\prox_f(x) = \arg\min_{z \in \Hm}\inv{2}\|z - x\|_2^2 + f(z)$ for
$f\in\Gamma_0(\Hm)$. If $f=\imath_C$ for some closed convex set
$C\subset \Hm$, $\prox_{f}(x)$ is equivalent to orthogonal projection
onto $C$. For $f(x)=\|x\|_1$, $\prox_{\gamma f}(x)$ is given by
component-wise soft-thresholding of $x$ by threshold $\gamma$
\cite{Fadili2009}. In addition, proximity operators of lsc convex
functions exhibit nice properties with respect to translation,
composition with frame operators, dilation,
etc. \cite{combettes2006srp,combettes2008pdm}.

In problem ({$\rm \bf P$}) with $f=f_1+f_2$, the resolvent operator
$\mathcal{J}_{\beta\partial{f}}=(\Id + \beta \partial f)^{-1}$
typically cannot be calculated in closed-form. Monotone operator
splitting methods do not attempt to evaluate this resolvent mapping
directly, but instead perform a sequence of calculations involving
separately the individual resolvent operators
$\mathcal{J}_{\beta\partial{f_1}}$ and
$\mathcal{J}_{\beta\partial{f_2}}$. The latter are hopefully easier to
evaluate, and this holds true for our functionals in BPDQ$_p$.

Since for BPDQ$_p$, both $f_1$ and $f_2$ are non-differentiable, we
use 
a particular monotone operator splitting method known as the
Douglas-Rachford (DR) splitting.
It can be written as the following compact recursion formula
\cite{Fadili2009}
\begin{equation}
\label{eq:DR-iter}
x^{(t+1)} = (1-\tfrac{\alpha_t}{2})\,x^{(t)} + 
\tfrac{\alpha_t}{2}\,S^\odot_{\gamma}\circ\Proj^\odot_{T_p(\epsilon)}(x^{(t)}),
\end{equation}
where $A^\odot \triangleq 2A - \Id$ for any operator $A$, $\alpha_t
\in (0,2)$ for all $t \in \mathbb{N}$, $S_{\gamma}=\prox_{\gamma f_1}$
is the component-wise  soft-thresholding operator with threshold
$\gamma>0$ and $\Proj_{T_p(\epsilon)} = \prox_{f_2}$ is the orthogonal
projection onto the tube $T^p(\epsilon)$. From
\cite{combettes2004smi}, one can show that the sequence $(x^{(t)})_{t
  \in \mathbb{N}}$ converges to some point $x^\star$ and
$\Proj_{T_p(\epsilon)}(x^\star)$ is a solution of BPDQ$_p$. In the
next Section, we provide a way to compute
$\Proj_{T_p(\epsilon)}(x^\star)$ efficiently.

\subsection{Proximity operator of the $\ell_p$ fidelity constraint}
Each step of the DR iteration \eqref{eq:DR-iter} requires computation
of $\prox_{f_2} = \Proj_{T^p(\epsilon)}$ for
$T^p(\epsilon)=\{x\in\Rbb^N:\|y_{\rm q} - \Phi x\|_p\leq \epsilon\}$.
We present an iterative method to compute this projection for 
%
$2\leq p\leq \infty$.

Notice first that, defining the unit $\ell_p$ ball $B^p
=\{y\in\Rbb^m:\norm{y}_p\leq 1\} \subset \Rbb^m$, we have 
$$
f_2(x)\ =\ \imath_{T^p(\epsilon)}(x) \ =\ (\imath_{B^p} \circ A_\epsilon)(x),
$$ 
with the affine operator $A_\epsilon(x) \triangleq \tinv{\epsilon}(\Phi x - y_{\rm q})$.

\begin{figure*}
  \centering
  \subfigure[\label{fig:first-exper-qual}]
  {\includegraphics[width=5.4cm]{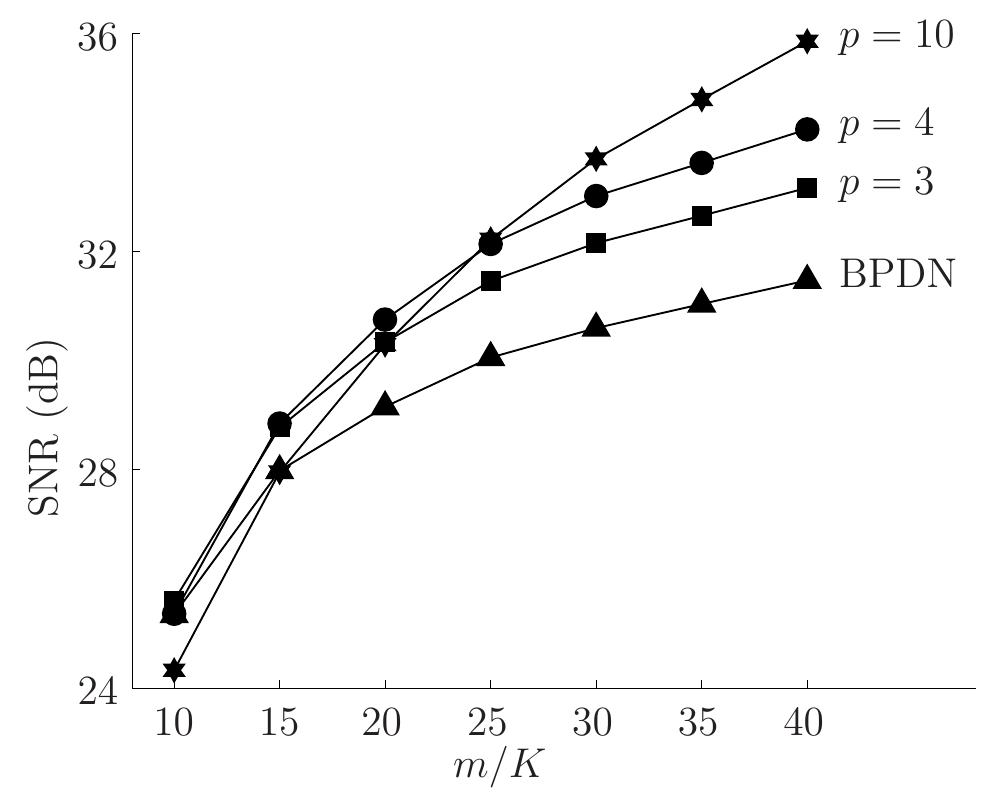}}
  \subfigure[\label{fig:first-exper-qual-std}]
  {\includegraphics[width=5.4cm]{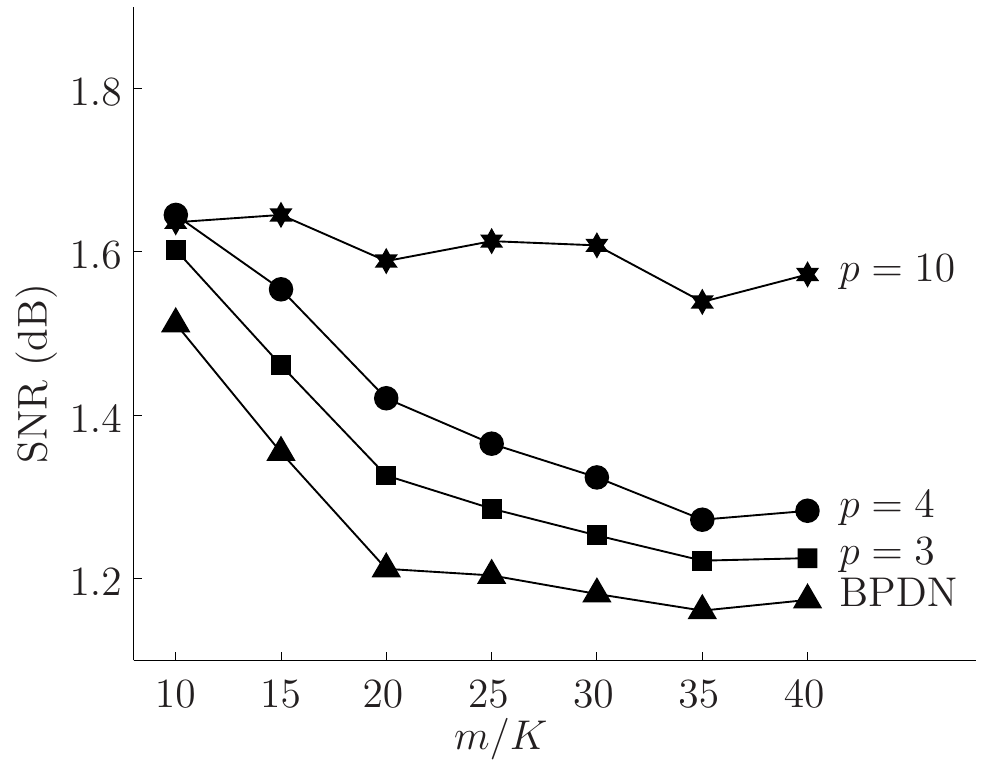}}  
  \subfigure[\label{fig:first-exper-qcfrac}]
  {\includegraphics[width=5.4cm]{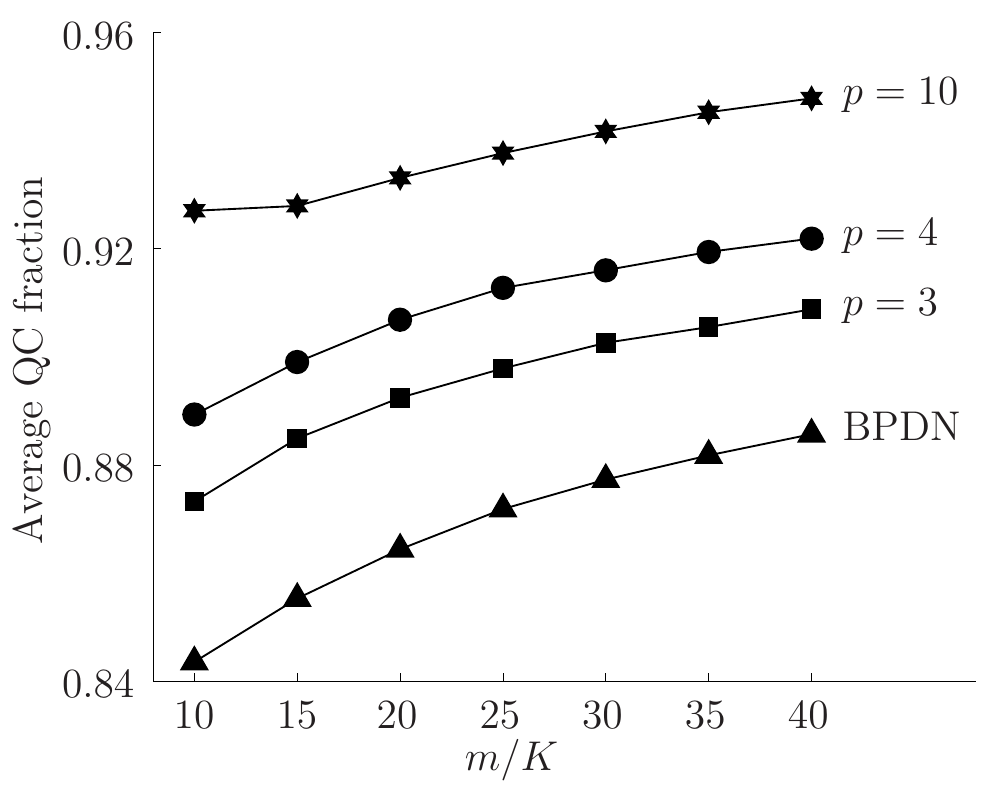}}  
\caption{Quality of BPDQ$_p$ for different $m/K$ and $p$. Mean (a)
  and standard deviation (b) of SNR. (c) Fraction of
  coefficients satisfying QC.}
  \label{fig:first-exper}
\end{figure*}

The proximity operator of a pre-composition of a function $f\in
\Gamma_0(\Hm)$ with an affine operator can be computed from the
proximity operator of $f$. Indeed, let $\Phi'\in\Rbb^{m\times N}$ and
the affine operator $A(x)\triangleq \Phi'x - y$ with $y\in\Rbb^m$. If
$\Phi'$ is a tight frame of $\Hm$, \ie $\Phi'\Phi'^* = c\Id$ for some
$c>0$, we have
\[
\prox_{f\circ A}(x) = x + c^{-1}\Phi'^* \big(\prox_{c f} - \Id\big)(A(x)) ~,
\]
\cite{combettes2006srp,Fadili2009}. Moreover, for a general bounded matrix $\Phi'$,
we can use the following lemma.
\begin{lemma}[\cite{Fadili2009}]
\label{lem:proxA}
Let $\Phi'\in\Rbb^{m\times N}$ be a matrix with bounds $0\leq c_1<c_2<\infty$ such that
$c_1\Id\leq \Phi'\Phi'^*\leq c_2\Id$ and let $\{\beta_t\}_{t\in\Nbb}$ be a
sequence with $0 < \inf_t \beta_t \leq \sup_t \beta_t < 2/c_2$. Define
\begin{equation}
\label{eq:proxFBdual}
\begin{array}{@{}r@{\,}c@{\,}l@{}}
u^{(t+1)}&=&\beta_t(\Id -
\prox_{\beta_t^{-1}f})
\big(\beta_t^{-1}u^{(t)} + A(p^{(t)})\big),\\[1mm]
p^{(t+1)} &=& x - \Phi'^* u^{(t+1)}.
\end{array}
\end{equation}
If the matrix $\Phi'$ is a general frame of $\mathcal{H}$,
\ie $0<c_1<c_2<\infty$, then $f \circ A \in \Gamma_{0}(\mathcal{H})$.
In addition, $u^{(t)} \rightarrow \bar{u}\in\Rbb^m$ and $p^{(t)}
\rightarrow \prox_{f\circ A}(x)=x - \Phi'^* \bar{u}$ in
\eqref{eq:proxFBdual}. More precisely, both $u^{(t)}$ and $p^{(t)}$
converge linearly and the best convergence rate is attained for
$\beta_t\equiv 2/(c_1+c_2)$ with $\norm{u^{(t)} - \bar{u}} \leq
\big(\tfrac{c_2-c_1}{c_2+c_1}\big)^t\norm{u^{(0)} - \bar{u}}$.
Otherwise, if $\Phi'$ is just bounded (\ie $c_1=0 < c_2 < \infty$),
and if $f\circ A \in \Gamma_{0}(\mathcal{H})$, apply
\eqref{eq:proxFBdual}, and then $u^{(t)} \rightarrow \bar{u}$ and
$p^{(t)} \rightarrow \prox_{f\circ A}(x)=x - \Phi'^* \bar{u}$ at the
rate $O(1/t)$.
\end{lemma}

In conclusion, computing $\prox_{f_2}$ may be reduced to applying the
orthogonal projection $\prox_{\imath_{B^p}} = \mathcal{P}_{B^p}$ by
setting $f=\imath_{B^p}$, $\Phi' = \Phi/\epsilon$ and $y=y_{\rm
  q}/\epsilon$ inside the iterative method (\ref{eq:proxFBdual}) with
a number of iterations depending on the selected application (see
Section \ref{sec:experiment}).

For $p=2$ and $p=\infty$, the projector $\mathcal{P}_{B^p}$ has an
explicit form.  Indeed, if $y$ is outside the closed unit
$\ell_p$-ball in $\Rbb^m$, then $\Proj_{B^2} (y) =
\frac{y}{\norm{y}_2}$; and $(\Proj_{B^\infty} (y))_i=\sign(y_i) \times \min\{1,|y_i|\}$ 
for $1 \leq i \leq m$.

Unfortunately, for $2<p<\infty$ no known closed-form for the
projection exists. Instead, we describe an iterative method. Set
$f_y(u) = \frac{1}{2}\norm{u-y}_2^2$ and $g(u)=\norm{u}_p^p$.  

If $\norm{y}_p \leq 1$, $\mathcal{P}_{B^p} (y) = y$. For $\norm{y}_p>1$,
the projection $\mathcal{P}_{B^p}$ is the solution of the constrained
minimization problem $u^\star = \arg\min_u f_y(u)\ \ {\rm s.t.}\ \ g(u) =
1$. Let $L(u,\lambda)$ be its Lagrange function (for $\lambda \in \Rbb$)
\begin{equation}
\label{eq-lagrange}
L(u,\lambda)\ =\ f_y(u) + \lambda\,(g(u) - 1).
\end{equation}
Without loss of generality, by symmetry, we may work in the positive\footnote
{The general solution can be obtained by appropriate axis mirroring.}
orthant $u_i\geq 0$ and $y_i\geq 0$, since the point $y$ and its
projection $u^\star$ belong to the same orthant of $\Rbb^m$,
\ie $y_iu^\star_i\geq 0$ for all $1\leq i\leq m$.

As $f_y$ and $g$ are continuously differentiable, the
Karush-Kuhn-Tucker system corresponding to \eqref{eq-lagrange} is
\begin{eqnarray}
\begin{split}
\label{eq-kkt}
\nabla_u L(u^\star,\lambda^\star) &= \nabla_u f_y(u^\star) + \lambda^\star \nabla_u g(u^\star) = 0 \\
\nabla_\lambda L(u^\star,\lambda^\star) &= g(u^\star) - 1 = 0 ~,
\end{split}
\end{eqnarray}
where the solution $u^\star$ is non-degenerate by strict convexity in
$u$ \cite{Nemirovski}, and $\lambda^\star$ the corresponding Lagrange
multiplier.

Let us write $z=(u,z_{m+1}=\lambda) \in\Rbb^{m+1}$ and
$F = \nabla_z L :\Rbb^{m+1} \to \Rbb^{m+1}$ as
\begin{equation*}
F_i(z)\ =\ \begin{cases}
\,z_i\ +\ p\,z_{m+1}\,z_i^{p-1}\ -\ y_i&\mbox{if}\ i\leq m,\\
\big( \sum_{j=1}^m z_j^p \big) -1&\mbox{if}\ i=m+1.
\end{cases}
\end{equation*}
The KKT system \eqref{eq-kkt} is equivalent to $F(z^\star)=0$, where
the desired projection $u^\star$ is then given by the first $m$
coordinates of $z^\star$. This defines a system of $m+1$ equations
with $m+1$ unknowns $(u^\star,\lambda^\star)$ that we can solved
efficiently with the Newton method. This is the main strategy
underlying sequential quadratic programming used to solve general-type
constrained optimization problems \cite{Nemirovski}.

Given an initialization point $z^0$, the successive iterates are defined by
\begin{equation} 
\label{eq-newtons-method}
z^{n+1} = z^n - V(z^n)^{-1} F(z^n),
\end{equation}
where $V_{ij} = \frac{\partial F_i}{\partial z_j}$ is the Jacobian
associated to $F$. If the iterates sequence $(z^{n})_{n\geq 0}$ is
close enough to $(u^\star,\lambda^\star)$, we known that the Jacobian
is nonsingular as $u^\star$ is non-degenerate. Moreover, since that
Jacobian has a simple block-invertible form, we may compute
(\cite{zwillinger2003csm}, p.125)
\begin{equation}
\label{eq-practical-jacobian-application}
V^{-1}(z)\ =\ \tinv{\mu} 
\begin{pmatrix}
\mu D^{-1}u\ +\ \big(z_{m+1} - \bar b^Tu\big)
\bar b\,\\
(\bar b^Tu - z_{m+1})
\end{pmatrix}\!,
\end{equation}
where $D\in \Rbb^{m\times m}$ is a diagonal matrix with $D_{ii}(z)=1\,
+\, p(p-1)z_{m+1}z_i^{p-2}$, $b\in \Rbb^m$ with $b_i(z)=pz_i^{p-1}$
for $1\leq i\leq m$, $\bar b = D^{-1}b$, $\mu= b^T D^{-1} b = \bar b^T
D\,\bar b$. This last expression can be computed efficiently as $D$ is
diagonal.

We initialize the first $m$ components of $z^0$ by the direct radial
projection of $y$ on the unit $\ell_p$-ball, $u^{0} = y/\norm{y}_p$,
and initialize $z^0_{m+1} = \arg\min_\lambda \|F(u^{0},\lambda)\|_2$.

In summary, to compute $\mathcal{P}_{B^p}$, we run
\eqref{eq-newtons-method} using
\eqref{eq-practical-jacobian-application} to calculate each update
step. We terminate the iteration when the norm of $\norm{F(z^n)}_2$
falls below a specified tolerance. Since the Newton method converges
superlinearly, we obtain error comparable to machine precision with
typically fewer than 10 iterations.

\section{Experiments}
\label{sec:experiment}

As an experimental validation of the BPDQ$_p$ method, we ran two sets of
numerical simulations for reconstructing signals from quantized
measurements. For the first experiment we studied recovery of exactly
sparse random 1-D signals, following very closely our theoretical
developments. Setting the dimension $N=1024$ and the sparsity level
$K=16$, we generated 500 $K$-sparse signals where the non-zero
elements were drawn from the standard Gaussian distribution $\mathcal{N}(0,1)$,
and located at supports drawn uniformly in $\{1,\,\cdots,N\}$.  For
each sparse signal $x$, $m$ quantized measurements were recorded as in
model (\ref{eq:quantiz-fmwk}) with a SGR matrix $\Phi\in\Rbb^{m\times
N}$. The bin width was set to $\upalpha = \|\Phi x\|_\infty/40$.

The decoding was accomplished with BPDQ$_p$ for various moments $p\geq
2$ using the optimization algorithm described in Section \ref{sec:implementations}.
In particular, the overall Douglas-Rachford procedure
\eqref{eq:DR-iter} was run for 500 iterations.  At each DR step, the
method in \eqref{eq:proxFBdual} was iterated until the relative error
$\tfrac{\norm{p^{(t)}-p^{(t-1)}}_2}{\norm{p^{(t)}}_2}$ fell below $10^{-6}$; the
required number of iterations was dependent on $m$ but was fewer than
700 in all cases examined.

In Figure~\ref{fig:first-exper}, we plot the average quality of the
reconstructions of BPDQ$_p$ for various values of $p\geq 2$ and
$m/K\in[10,40]$.  We use the quality measure $\mbox{SNR}(\hat{x};x) =
20\log_{10}\frac{\|x\|_2}{\|x-\hat{x}\|_2}$, where $x$ is the true
original signal and $\hat{x}$ the reconstruction.  As can be noticed,
at higher oversampling factors $m/K$ the decoders with higher $p$ give
better reconstruction performance. Equivalently, it can also be
observed that at lower oversampling factors, increasing $p$ beyond a
certain point degrades the reconstruction performance. These two
effects are consistent with the remarks noted at the end of Section
\ref{sec:BPDQ-approx-error-for-quantiz}, as the sensing matrices may
fail to satisfy the RIP$_{p,2}$ if $p$ is too large for a given
oversampling factor.

One of the original motivations for the BPDQ$_p$ decoders is that they are
closer to enforcing quantization consistency than the BPDN decoder. To
check this, we have examined the ``quantization consistency
fraction'', \ie the average fraction of remeasured coefficients
$(\Phi \hat{x})_i$ that satisfy $|(\Phi \hat{x})_i-y_i| <
\tfrac{\upalpha}{2}$. These are shown in Figure \ref{fig:first-exper}
(c) for various $p$ and $m/K$. As expected, it can be clearly seen that increasing
$p$ increases the QC fraction.

\begin{figure}
\centering
  {\includegraphics[height=3.5cm]{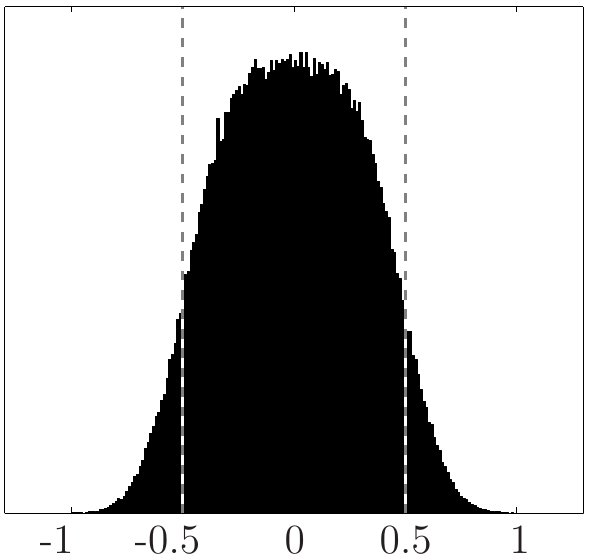}}\hspace{.2cm} 
  {\includegraphics[height=3.5cm]{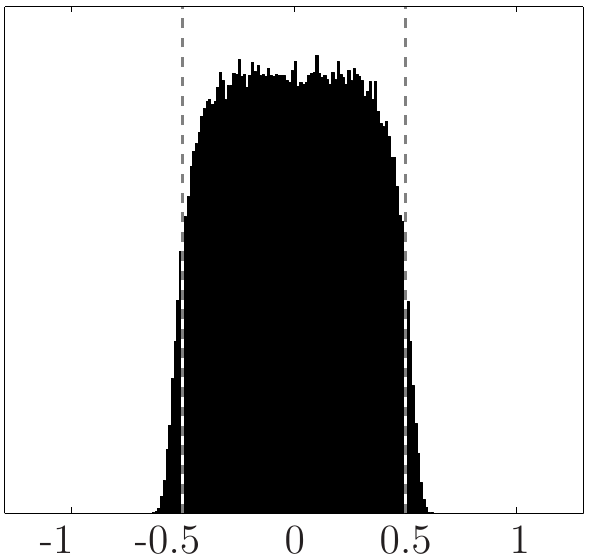}}  
\caption{Histograms of $\upalpha^{-1}(\Phi \hat{x}-y)_i$. Left, $p=2$.
Right, $p=10$.
\label{fig:normalized-residual-histograms} }       
\end{figure}

An even more explicit illustration of this effect is afforded by
examining histograms of the normalized residual $\upalpha^{-1}(\Phi
\hat{x}-y)_i$ for different $p$. For reconstruction exactly satisfying
QC, these normalized residuals should be supported on $[-1/2, 1/2]$. In
Figure \ref{fig:normalized-residual-histograms} we show histograms
of normalized residuals for $p=2$ and $p=10$, for the case
$m/K=40$. The histogram for $p=10$ is indeed closer
to uniform on $[-1/2, 1/2]$.

For the second experiment, we apply a modified version of the BPDQ$_p$ to
an undersampled MRI reconstruction problem. Using an example
similar to \cite{lustig2007sma}, the original image is a $256\times256$
pixel ``synthetic angiogram'', \ie $N=256^2$, comprised of 10
randomly placed non-overlapping ellipses. The linear measurements are
the real and imaginary parts of a fraction $\rho$ of the Fourier coefficients
at randomly selected locations in Fourier space, giving $m = \rho N$
independent measurements. These random locations form the index set
$\Omega\subset \{1,\,\cdots,N\}$ with $|\Omega|=m$. Experiments were
carried out with $\rho\in\{1/6,1/8,1/12\}$, but we show results only
for $\rho = 1/8$. These
were quantized with a bin width $\upalpha=50$, giving at most 12
quantization levels for each measurement.

For this example, we  modify the BPDQ$_p$ program \ref{eq:BPDQ} by replacing the
$\ell_1$ term by the total variation (TV) semi-norm \cite{Rudin1992}.
This yields the problem\\[-2mm]
\begin{equation*}
\argmin_u \|u\|_{TV}\ \ {\rm s.t.}\ \ \|y - \Phi u\|_p \leq \epsilon,
\vspace{-2mm}
\end{equation*}
where $\Phi = F_{\Omega}$ is the restriction of Discrete Fourier
Transform matrix $F$ to the rows indexed in $\Omega$. 

This may be solved with the Douglas-Rachford iteration
\eqref{eq:DR-iter}, with the modification that $S_\gamma$ be replaced
by the proximity operator associated to $\gamma$ times the TV norm,
\ie by $\prox_{\gamma\|\cdot\|_{\mathrm{TV}}}(y) =
\argmin_u\frac{1}{2}\|y-u\|^2 + \gamma\|u\|_{\mathrm{TV}}$. The latter
is known as the Rudin-Osher-Fatemi model, and numerous methods exist
for solving it exactly, including
\cite{Vogel1996,Chan1999,Chambolle2004,Wang2008}. In this work, we use
an efficient projected gradient descent algorithm on the dual problem,
see \eg \cite{Fadili2009}. Note that the sensing matrix $F_\Omega$ is
actually a tight frame, \ie $F_\Omega F_\Omega^*=\Id$, so we do not
need the nested inner iteration \eqref{eq:proxFBdual}.

\begin{figure}
  \centering
  \includegraphics[height=4.05cm]{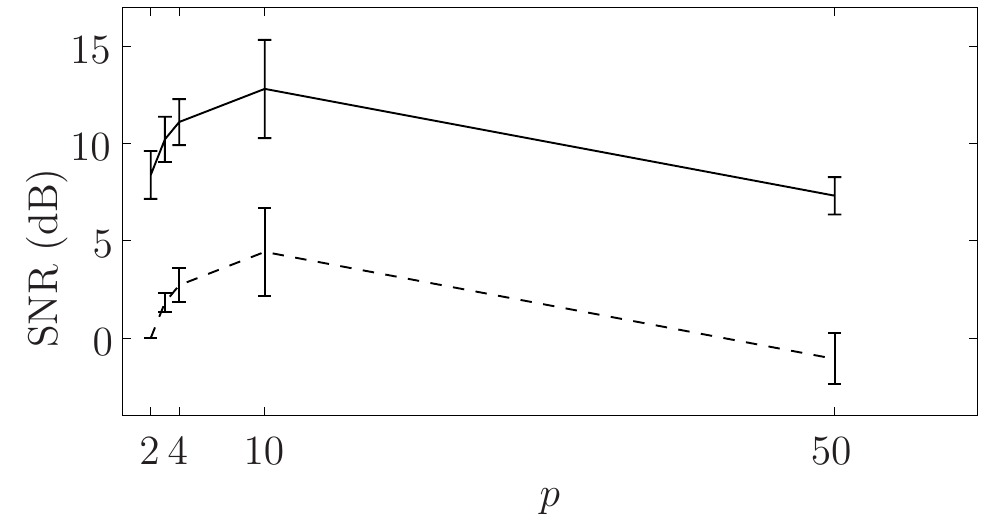}
  \caption{\label{fig:2d-snr}Average SNR (solid) and SNR improvement over
    BPDN (dashed) as a function of $p$, for the synthetic
    angiogram reconstruction simulations. Error bars indicate 1
    standard deviation. }
\end{figure}

\begin{figure*}
  \centering
  \subfigure[\label{fig:second-exper-orig}Original image]
  {\includegraphics[height=5cm]{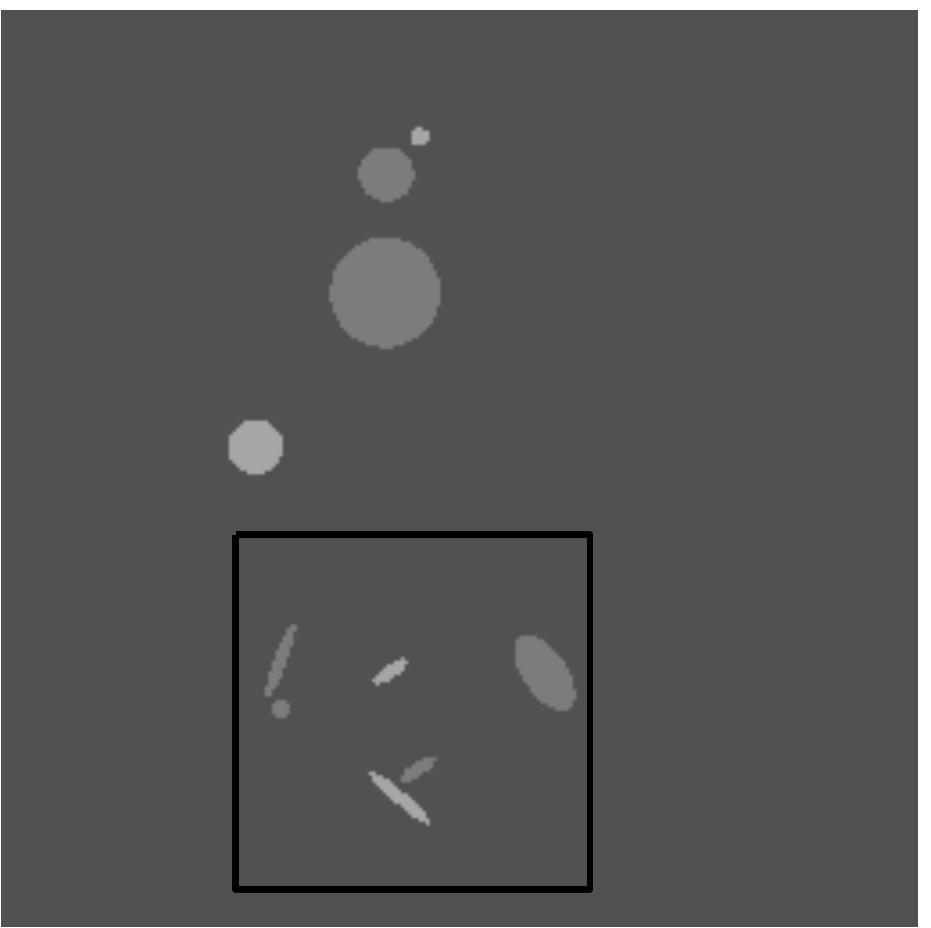}}\hspace{5mm} 
  \subfigure[\label{fig:second-exper-bpdn}SNR = 8.96 dB]
  {\includegraphics[height=5cm]{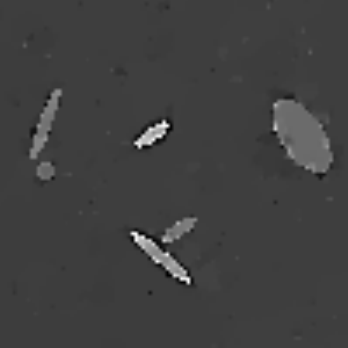}}\hspace{5mm} 
  \subfigure[\label{fig:second-exper-bpdq}SNR = 12.03 dB]
  {\includegraphics[height=5cm]{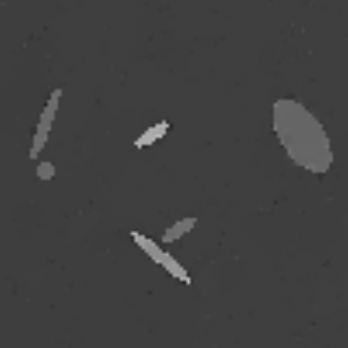}}
  \caption{\label{fig:examp2d}Reconstruction of synthetic angiograms
    from undersampled Fourier measurements, using TV
    regularization. (a) Original, showing zoom area (b) BPDN (zoom)
    (c) BPDQ$_{10}$ (zoom).}
\end{figure*}

We show the SNR of the BPDQ$_p$ reconstructions as a function of $p$
in Figure~\ref{fig:2d-snr}, averaged over 50 trials where both the
synthetic angiogram image and the Fourier measurement locations are
randomized. This figure also depicts the SNR improvement of
BPDQ$_p$-based reconstruction over BPDN.  For these simulations we
used 500 iterations of the Douglas-Rachford recursion
\eqref{eq:DR-iter}. This quantitative results are confirmed by
  visual inspection of Figure~\ref{fig:examp2d}, where we compare
  $100 \times 100$ pixel details of the reconstruction results
  with BPDN and with BPDQ$_p$ for $p=10$, for one particular
  instance of the synthetic angiogram signal.

  Note that this experiment lies far outside of the justification
  provided by our theoretical developments, as we do not have any
  proof that the sensing matrix $F_\Omega$ satisfies the RIP$_{p,2}$,
  and our theory was developed only for $\ell_1$ synthesis-type
  regularization, while the TV regularization is of analysis type.
  Nonetheless, we obtain results analogous to the previous 1-D
  example; the BPDQ$_p$ reconstruction shows improvements both in SNR
  and visual quality compared to BPDN. These empirical results suggest
  that the BPDQ$_p$ method may be useful for a wider range of
  quantized reconstruction problems, and also provoke interest for
  further theoretical study.

\section{Conclusion and Further Work}
\label{sec:conclusion}

The objective of this paper was to show that the BPDN reconstruction
program commonly used in Compressed Sensing with noisy measurements is
not always adapted to quantization distortion. We introduced a new
class of decoders, the Basis Pursuit DeQuantizers, and we have shown
both theoretically and experimentally that BPDQ$_p$ exhibit a
substantial reduction of the reconstruction error in oversampled
situations. 

A first interesting question for further study would be to
characterize the evolution of the optimal moment $p$ with the
oversampling ratio. This would allow for instance the selection of the
best BPDQ$_p$ decoder in function of the precise CS coding/decoding
scenario. Second, it is also worth investigating the existence of
other RIP$_{p,2}$ random matrix constructions, \eg using the Random
Fourier Ensemble. Third, a more realistic coding/decoding scenario
should set $\upalpha$ theoretically in function of the bit budget
(rate) available to quantize the measurements, of the sensing matrix
and of some a priori on the signal energy. This should be linked also
to the way our approach can integrate the saturation of the quantized
measurements \cite{laska09}. Finally, we would like to extend our
approach to non-uniform scalar quantization of random measurements,
generalizing the quantization consistency and the optimization
fidelity term to this more general setting.

\section{Acknowledgments}
\label{sec:thanks}

LJ and DKH are very grateful to Prof. Pierre Vandergheynst (Signal
Processing Laboratory, LTS2/EPFL, Switzerland) for his useful advices
and his hospitality during their postdoctoral stay in EPFL. 

\appendices

\section{Proof of Proposition \ref{prop:gauss-rip-inf}}
\label{app:why-rip-p}

Before proving Proposition \ref{prop:gauss-rip-inf}, let us recall
some facts of measure concentrations \cite{ledoux1991pbs,
  ledoux2001cmp}.

In particular, we are going to use the concentration property of any
Lipschitz function over $\Rbb^m$, \ie $F$ such that $\|F\|_{\rm Lip}
\triangleq \sup_{u,v\,\in\,\Rbb^m,\ u \neq v}\,\frac{|F(u) -
  F(v)|}{\|u-v\|_2} < \infty$.  If $\|F\|_{\rm Lip}\leq 1$, $F$ is
said 1-Lipschitz.

\begin{lemma}[Ledoux, Talagrand \cite{ledoux1991pbs} (Eq. 1.6)] 
\label{lem:measconc-lipschitz}
If $F$ is
Lipschitz with $\lambda=\|F\|_{\rm Lip}$, then, for the random vector $\xi\in\Rbb^m$ with $\xi_i\sim_{\rm iid} \mathcal{N}(0,1)$, 
\begin{equation*}
P_\xi\big[\,| F(\xi) - \mu_F | > r\,\big]\ \leq\ 2
e^{-\inv{2} r^2 \lambda^{-2}},\quad {\rm for}\ r>0,
\end{equation*}
with $\mu_F = \E F(\xi) = \int_{\Rbb^m} F(x)\, \gamma^m(x)\,\ud^m x$ and
$\gamma^m(x) = (2\pi)^{-m/2}\ e^{-\|x\|_2^2/2}$.
\end{lemma}

A useful tool that we will use is the concept of a \emph{net}. An
$\epsilon$-net ($\epsilon > 0$) of $A\subset\Rbb^K$ is a subset
$\mathcal{S}$ of $A$ such that for every $t\in A$, one can find
$s\in\mathcal{S}$ with $\|t-s\|_2\leq \epsilon$. In certain cases, the
size of a $\epsilon$-net can be bounded.

\begin{lemma}[\cite{ledoux2001cmp}] 
\label{lem:size-eps-net-sphere}
There exists a $\epsilon$-net $\mathcal{S}$ of the unit sphere of
$\Rbb^K$ of size $|\mathcal{S}|\leq (1+\frac{2}{\epsilon})^K$.
\end{lemma}

\noindent We will use also this fundamental result.
\begin{lemma}[\cite{ledoux2001cmp}] 
\label{lemma:eps-net-sufficient}
Let $\mathcal{S}$ be a $\epsilon$-net of the unit sphere in
$\Rbb^K$. Then, if for some vectors $v_1,\cdots,v_K$ in the Banach
space $B$ normed by $\|\!\cdot\!\|_B$, we have $1-\epsilon \leq
\big\|\sum_{i=1}^K s_iv_i \big\|_B\leq 1+\epsilon$ for all
$s=(s_1,\cdots,s_K)\in\mathcal{S}\subset\Rbb^K$, then
\begin{equation*}
(1 - \beta)\,\|t\|_2\ \leq \big\|\sum_{i=1}^K
t_iv_i \big\|_B\ \leq \ (1+\beta)\,\|t\|_2,
\end{equation*} 
for all $t\in\Rbb^K$, with $\beta = \tfrac{2\epsilon}{1-\epsilon}$.
\end{lemma}
\noindent In our case, the Banach space $B$ is $\ell_p^m=(\Rbb^m,\|\!\cdot\!\|_p)$ for
$1\leq p\leq \infty$, i.e~$\Rbb^m$ equipped with the norm
$\|u\|_p^p=\sum_i |u_i|^p$. With all these concepts, we can now
demonstrate the main proposition.
\begin{proof}[Proof of Proposition \ref{prop:gauss-rip-inf}]
  Let $p\geq 1$.  We must prove that for a SGR matrix
  $\Phi\in\Rbb^{m\times N}$, \ie with $\Phi_{ij}\sim_{iid} \mathcal{N}(0,1)$,
  with the right number of measurements $m$, there exist a radius
  $0<\delta<1$ and a constant $\mu_{p,2}>0$ such that
  \begin{equation}
    \label{eq:rip-p-annex}
    \mu_{p,2}\,\sqrt{1-\delta}\,\|x\|_2\ \leq\ \|\Phi x\|_p\ \leq\
    \mu_{p,2}\,\sqrt{1+\delta}\,\|x\|_2,
  \end{equation}
  for all $x\in\Rbb^N$ with $\|x\|_0\leq K$.
  
  We begin with a unit sphere $S_T=\{u\in\Rbb^N: {\rm supp}\,u=T,
  \|u\|_2=1\}$ for a fixed support $T\subset\{1,\cdots,N\}$ of size
  $|T|=K$. Let $\mathcal{S}_T$ be an $\epsilon$-net of $S_T$.  We
  consider the SGR random process that generates $\Phi$ and, by an
  abuse of notation, we identify it for a while with $\Phi$ itself. In
  other words, we define the random matrix $\Phi =
  (\Phi_1,\cdots,\Phi_N)\in\Rbb^{m\times N}$ where, for all $1\leq
  i\leq N$, $\Phi_j\in\Rbb^m$ is a random vector of probability
  density function (or \emph{pdf}) $\gamma^m(u)=\Pi_{i=1}^m
  \gamma(u_i)$ for $u\in\Rbb^m$ and $\gamma(u_i) = \inv{\sqrt{2\pi}}
  e^{-u_i^2/2}$ (the standard Gaussian pdf). Therefore, $\Phi$ is
  related to the pdf $\gamma_\Phi(\phi)=\Pi_{j=1}^N \gamma^m(\phi_j)$,
  $\phi=(\phi_1,\cdots,\phi_N)\in\Rbb^{m\times N}$. 

  Since the Frobenius norm $\|\phi\|_{\mathcal
    F}=(\sum_{jk}|\phi_{jk}|^2)^{1/2}$ of $\phi$ and the pdf
  $\gamma_\Phi(\phi)\propto e^{-\|\phi\|^2_{\mathcal F}/2}$ are invariant
  under a global rotation in $\Rbb^N$ of all the rows of $\phi$, it is
  easy to show that for unit vector $s\in\Rbb^N$, $ P_\Phi\big[|
  F(\Phi s) - \mu_F | > r\big] = P_\Phi\big[| F(\Phi_1) - \mu_F | >
  r\big] \leq 2 e^{-\inv{2} r^2 \lambda^{-2}} $, using Lemma
  \ref{lem:measconc-lipschitz} on the SGR vector $\Phi_1$.

The above holds for a single $s$. To obtain a result valid for all
$s\in \mathcal{S}_T$ we may use the union bound. As
$|\mathcal{S}_T|\leq (1 + 2/\epsilon)^K$ by Lemma
\ref{lem:size-eps-net-sphere}, setting $r=\epsilon\mu_F$ for
$\epsilon>0$, we obtain
\begin{equation*}
P_\Phi\big[\,|\,\mu_F^{-1}\,F(\Phi s) - 1 | > \epsilon\big]\ \leq\ 2
\,e^{K\log(1+2\epsilon^{-1})-\inv{2} \epsilon^2\mu_F^2 \lambda^{-2}},
\end{equation*}
for all $s\in \mathcal{S}_T$.

Taking now $F(\cdot)=\|\!\cdot\!\|_p$ for $1\leq p\leq \infty$, we
have $\mu_F=\mu_{p,2}=\E\|\xi\|_p$ for a SGR vector $\xi\in\Rbb^m$. The
Lipschitz value is $\lambda=\lambda_p=1$ for $p\geq 2$, and $\lambda =
\lambda_p = m^{\frac{2-p}{2p}}$ for $1\leq p\leq 2$. Consequently,
\begin{equation}
\label{lemma:conc-bound-enet}
(1-\epsilon)\ \leq\ \norm{ \tfrac{1}{\mu_{p,2}} \Phi s}_p\ \leq (1 + \epsilon), 
\end{equation}
for all $s\in \mathcal{S}_T$, with a probability higher than $1 - 2
\,\exp(K\log(1+2\epsilon^{-1}) -\tinv{2}\epsilon^2\mu_{p,2}^2\lambda_p^{-2})$. 

We apply Lemma \ref{lemma:eps-net-sufficient} by noting that, as $s$
has support of size $K$, (\ref{lemma:conc-bound-enet}) may be written
as
\begin{equation*}
1-\epsilon\ \leq\ \norm{ \sum_{i=1}^K s_i v_i }_p\ \leq\ 1+\epsilon
\end{equation*}
where $v_i\in\Rbb^m$ are the columns of $\frac{1}{\mu_{p,2}} \Phi$ corresponding to
the support of $s$ (we abuse notation to let $s_i$ range only over the
support of $s$). Then according to Lemma
\ref{lemma:eps-net-sufficient} we have, with the same probability
bound and for $(\sqrt{2}-1)\delta =
\tfrac{2\epsilon}{1-\epsilon}$,
\begin{multline} 
\label{lemma:l2bound-Ksupp}
\sqrt{1-\delta}\,\|x\|_2\ \leq\ (1-(\sqrt{2}-1)\delta)\,\|x\|_2\ \leq\
\|\Phi x\|_p\\ \leq\ (1 + (\sqrt{2}-1)\delta)\,\|x\|_2\ \leq\ \sqrt{1+\delta}\,\|x\|_2, 
\end{multline}
for all $x\in\Rbb^N$ with ${\rm supp}\, x = T$.

The result can be made independent of the choice of $T\subset
\{1,\,\cdots,N\}$ by considering that there are $\binom{N}{K} \leq
(eN/K)^K$ such possible supports. Therefore, applying again an union
bound, (\ref{lemma:l2bound-Ksupp}) holds for all $K$-sparse $x$ in
$\Rbb^N$ with a probability higher than $1 - 2
\,e^{-\inv{2}\epsilon^2\mu_{p,2}^2\lambda_p^{-2} + K\log [e\frac{N}{K}(1+2\epsilon^{-1})]}$.

Let us bound this probability first for $1\leq p <\infty$. For $m\geq
\beta^{-1}\,2^{p+1}$ and $\beta^{-1}=p-1$, Lemma
\ref{lem:strict-bounds-mu_p} (page \pageref{lem:strict-bounds-mu_p})
tells us that $\mu_{p,2} \geq \tfrac{p-1}{p} \nu_p\,m^{\frac{1}{p}}$
with
$\nu_p=\sqrt{2}\,\pi^{-\frac{1}{2p}}\,\Gamma(\tfrac{p+1}{2})^{\frac{1}{p}}$.
A probability of success $1-\eta$ with $\eta<1$ is then guaranteed if
we select, for $1\leq p < 2$,
\begin{equation*}
m > \tfrac{2}{\epsilon^2 \nu_p^2}(\tfrac{p}{p-1})^2\big(K \log[e\tfrac{N}{K}(1
+ 2\epsilon^{-1})] + \log \tfrac{2}{\eta}\big),
\end{equation*}
since $\lambda_p = m^{\frac{2-p}{2p}}$, and for $2\leq p < \infty$,
\begin{equation}
\label{eq:p-g-2-meas-req}
m^{\frac{2}{p}} > \tfrac{2}{\epsilon^2 \nu_p^2}(\tfrac{p}{p-1})^2\big(K \log[e\tfrac{N}{K}(1
+ 2\epsilon^{-1})] + \log \tfrac{2}{\eta}\big),
\end{equation}
since $\lambda_p=1$. 

From now, $A \geq c\,B$ or $A \leq c\,B$ means that there exists a
constant $c>0$ such that these inequalities hold. According to the
lower bound found in Section \ref{sec:BPDQ-approx-error-for-quantiz},
$\nu_p > c\,\sqrt{p+1}$ implying that $\nu_p^{-2} \leq c$. Since
$(p/(p-1))^2\leq 4$ for any $p\geq 2$ and $\epsilon^{-1}\leq
\frac{\sqrt 2 + 1}{\sqrt 2 - 1}\delta^{-1} \leq 6\,\delta^{-1}$, we
find the new sufficient conditions,
\begin{equation*}
  m > c\,\delta^{-2} (\tfrac{p}{p-1})^2\big(
  K \log[e\tfrac{N}{K}(1
  + 12\delta^{-1})] + \log \tfrac{2}{\eta}\big),
\end{equation*}
for $1\leq p < 2$, and
\begin{equation*}
m^{2/p} > c \,\delta^{-2}\,\big(K \log[e\tfrac{N}{K}(1
  + 12\delta^{-1})] + \log \tfrac{2}{\eta}\big),
\end{equation*}
for $2\leq p < \infty$.

Second, in the specific case where $p=\infty$, since there exists a
$\rho>0$ such that $\mu_{\infty,2}\,\geq \rho^{-1} \sqrt{\log m}$,
with $\lambda_\infty=1$,
$
\log m\ >\ c \,\delta^{-2}\,\big(K \log[e\tfrac{N}{K}(1
  + 12\delta^{-1})]\ +\ \log \tfrac{2}{\eta}\big).
$
\end{proof}

Let us make some remarks about the results and the requirements of the
last proposition. Notice first that for $p=2$, we find the classical
result proved in \cite{JLmeetCS}. Second, as for the comparison
between the common RIP$_{2,2}$ proof \cite{JLmeetCS} and the tight
bound found in \cite{DoTa09}, the requirements on the measurements
above are possibly pessimistic, \ie the exponent $2/p$ occurring in
(\ref{eq:p-g-2-meas-req}) is perhaps too small.  Proposition
\ref{prop:gauss-rip-inf} has however the merit to prove that random
Gaussian matrices satisfy the RIP$_{p,2}$ in a certain range of
dimensionality.

\section{Link between $\delta$ and $m$ for SGR RIP$_{p,2}$ matrices}
\label{sec:delta-propto-m1p}

For $2\leq p< \infty$, Proposition \ref{prop:gauss-rip-inf} shows
that, if $\delta^{2} \geq c \,m^{-2/p}\,\big(K \log[e\tfrac{N}{K}(1 +
12\delta^{-1})] + \log \tfrac{2}{\eta}\big)$ for a certain constant
$c>0$, a SGR matrix $\Phi\in\Rbb^{m\times N}$ is RIP$_{p,2}$ of order
$K$ and radius $0<\delta<1$ with a probability higher than $1-\eta$.
Let us assume that $\delta > d m^{-1/p}$ for some $d>0$. We have,
$\log\delta^{-1} < \inv{p}\log m - \log d$, and therefore, the same
event occurs with the same probability bound when $\delta^{2} \geq c
\,m^{-2/p}\,\big(K \log[13 e\tfrac{N}{K}] + \frac{K}{p}\log m - K\log
d + \log \tfrac{2}{\eta}\big)$.  For high $m$ and for fixed $K,N$ and
$\eta$, this provides $\delta = O(m^{-1/p}\sqrt{\log m})$, which meets
the previous assumption.

\section{}
\label{sec:proof-lemma-refl}

\begin{proof}[Proof of Lemma \ref{lem:strict-bounds-mu_p}]
  The result for $p=\infty$ is due to \cite{ledoux1991pbs} (see Eq
  (3.14)).  Let $\xi\in\Rbb^m$ be a SGR vector, \ie $\xi_i\sim_{iid}
  \mathcal{N}(0,1)$ for $1\leq i\leq m$, and $1\leq p <\infty$.  First, the inequality
  $\E\|\xi\|_p\leq (\E\|\xi\|_p^p)^{1/p}$ follows from the application
  of the Jensen inequality $\varphi(\E\|\xi\|_p) \leq
  \E\varphi(\|\xi\|_p)$ with the convex function
  $\varphi(\cdot)=(\cdot)^p$. Second, the lower bound on
  $\E\|\xi\|_p$ arises from the observation that for
  $f:\Rbb^+\to\Rbb^+$ with $f(t)=t^{\inv{p}}$, and for a given
  $t_0>0$,
\begin{equation}
\label{eq:mino-ineq-expec}
f(t) \geq f(t_0) + f'(t_0)(t-t_0) +
pf''(t_0)(t-t_0)^2,
\end{equation}
for all $t\geq 0$.

Indeed, observe first that since $f^{(n)}(\alpha t') =
\alpha^{\inv{p}-n} f^{(n)}(t')$ for $\alpha > 0$ and $n\in\Nbb$, it is
sufficient to prove the result for $t_0=1$. Proving
\eqref{eq:mino-ineq-expec} amounts then to prove
$f(t) = t^{\inv{p}} \geq \tfrac{2p-1}{p}t - \tfrac{p - 1}{p}t^2$,
or equivalently,
$t^{\inv{p}-1} + \tfrac{p - 1}{p} t \geq \tfrac{2p-1}{p}$.
The LHS of this last inequality takes its minimum in $t=1$ with value
$\frac{2p-1}{p}$, which provides the result.

Since $\mu_{p,2} = \E\|\xi\|_p=\E f(\|\xi\|^p_p)$ and $\E(\|\xi\|^p_p -
\bar\mu_{p,2}) = 0$, using \eqref{eq:mino-ineq-expec} we find 
\begin{equation*}
\mu_{p,2} \geq 
(t_0)^{\frac{1}{p}-2}\big((2-\tinv{p})\bar\mu_{p,2} t_0 + 
(\tinv{p}-1)(\bar\mu_{p,2}^2+\bar\sigma_p^2)\big)
\end{equation*}
writing $\bar \mu_{p,2} = \E\|\xi\|_p^p$ and $\bar\sigma_p^2 =
\E(\|\xi\|^p_p-\bar\mu_{p,2})^2 = {\rm Var}\|\xi\|^p_p$. The RHS of the
last inequality is maximum for $t_0 = \bar\mu_{p,2}\,(1 +
\bar\mu_{p,2}^{\,-2}\,\bar\sigma_p^2)$. For that value, we get finally
\begin{equation*}
  \mu_{p,2}\ \geq\ (\E\|\xi\|^p_p)^{\inv{p}}\ 
\big(1\ +\ (\E\|\xi\|_p^p)^{-2}\,{\rm Var}\|\xi\|_p^p\big)^{\frac{1}{p}-1}.
\end{equation*}

Because of the decorrelation of the components of $\xi$, the last
inequality simplifies into
\begin{equation*}
\mu_{p,2}\ \geq\ m^{\inv{p}}\,(\E|g|^p)^{\inv{p}}\,
\big(\,1\ +\ m^{-1}(\E|g|^p)^{-2}\,{\rm Var}|g|^p\,\big)^{\frac{1}{p}-1},
\end{equation*}
with $g\sim \mathcal{N}(0,1)$.

Moreover, since $\E|g|^p =
2^{\frac{p}{2}}\,\pi^{-\inv{2}}\,\Gamma(\tfrac{p+1}{2})$ and using the
following approximation of the Gamma function \cite{spira1971cgf}
$|\Gamma(x) - (\tfrac{2\pi}{x})^{\frac{1}{2}}(\tfrac{x}{e})^x|\leq
\tinv{9 x}\,(\tfrac{2\pi}{x})^{\frac{1}{2}}(\tfrac{x}{e})^x$, valid
for $x\geq 1$, we observe that
$$
0.9\,(\tfrac{2\pi}{x})^{\inv{2}}(\tfrac{x}{e})^x\ \leq\
\Gamma(x)\ \leq\ 1.1\,(\tfrac{2\pi}{x})^{\frac{1}{2}}(\tfrac{x}{e})^x,
$$
that holds also if $x=\tfrac{p+1}{2}$ with $p\geq 1$.  Therefore,
$(\E |g|^p)^{-2}\,{\rm Var}|g|^p \leq
\big(\tfrac{1.1}{0.9^2}(\tfrac{e}{2})^{\inv{p}}(\tfrac{2p+1}{p+1})^{p} - 1\big) 
\leq
\tfrac{1.1}{0.9^2}(\tfrac{e}{2})^{\inv{2}}2^{p}$
and finally
\begin{equation*}
\mu_{p,2}\ \geq\ m^{\inv{p}}\,(\E|g|^p)^{\inv{p}}\,
\big(1\ +\ c\,\tfrac{2^p}{m}\,\big)^{\frac{1}{p}-1}
\end{equation*}
for a constant $c = \tfrac{1.1}{0.9^2}\,(\tfrac{e}{2})^{\inv{2}}<1.584<2$ 
independent of $p$ and $m$.
\end{proof}

\section{}
\label{sec:proof-lemma-bound-scp-lp}

\begin{proof}[Proof of Lemma \ref{lemma:bound-scp-lp}]
Notice first that since $J(\lambda w)=\lambda\,J(w)$ for any
$w\in\Rbb^m$ and $\lambda\in\Rbb$, it is sufficient to prove the
result for $\|u\|_2=\|v\|_2=1$.

The Lemma relies mainly on the geometrical properties of the Banach
space $\ell_p^m=(\Rbb^m,\|\!\cdot\!\|_p)$ for $p\geq 2$. In
\cite{bynum1976wpl,xu1991ibs}, it is explained that this space is
$p$-convex and 2-smooth. The smoothness involves in particular
\begin{equation}
\label{eq:smooth-ineq}
\|x+y\|_p^2 \leq \|x\|_p^2 + 2\scp{J(x)}{y} +
 (p-1)\|y\|_p^2,
\end{equation}
where $J=J_2$ and $J_r$ is the \emph{duality} mapping of \emph{gauge
  function} $t\to t^{r-1}$ for $r\geq 1$. For the Hilbert space
$\ell_2$, the relation \eqref{eq:smooth-ineq} reduces of course to the
\emph{polarization identity}. For $\ell_p$, $J_r$ is the differential
of $\tinv{r}\|\!\cdot\!\|_p^r$, \ie $(J_r(u))_i =
\|u\|^{r-p}\,|u_i|^{p-1}\,\sign u_i$. 

The smoothness inequality \eqref{eq:smooth-ineq} involves
\begin{equation}
\label{eq:scp_xy_first}
2\,\scp{J(x)}{y}\ \leq\ \|x\|_p^2\ +\ (p-1)\,\|y\|_p^2 - \|x-y\|_p^2,
\end{equation}
where we used the change of variable $y\to -y$.  

Let us take $x=\Phi u$ and $y=t\Phi v$ with $\|u\|_0=s$,
$\|v\|_0=s'$, $\|u\|_2=\|v\|_2=1$, $\supp u\, \cap\, \supp v = \emptyset$
and for a certain $t>0$ that we will set later. Because $\Phi$ is
assumed RIP$_{p,2}$ for $s$, $s'$ and $s+s'$ sparse signals, we deduce
\begin{multline*}
2\,\mu_{p,2}^{-2}\ t\ |\scp{J(\Phi u)}{\Phi v}|\ \leq\ 
(1+\delta_{s})\ + \\
(p-1)(1+\delta_{s'})t^2\ -\ (1-\delta_{s+s'})(1+t^2),
\end{multline*}
where the absolute value on the inner product arises from the
invariance of the RIP bound on \eqref{eq:scp_xy_first} under the
change $y\to -y$.  The value $\mu_{p,2}^{-2}|\scp{J(\Phi u)}{\Phi v}|$ is
thus bounded by an expression of type
$f(t)=\frac{\alpha+\beta t^2}{t}$ with $\alpha,\beta>0$ for
$p\geq 2$ given by $\alpha=\delta_{s}+\delta_{s+s'}$ and
$\beta=(p-2)+(p-1)\delta_{s'}+\delta_{s+s'}$.  Since the minimum of
$f$ is $2\sqrt{\alpha\beta}$, we get
\begin{multline}
\label{eq-pl3-a}
\mu_{p,2}^{-2}\,|\scp{J(\Phi u)}{\Phi v}|\ \leq\\
\big[(\delta_{s}+\delta_{s+s'})\big(\bar p +
\bar p\,\delta_{s'} + \delta_{s'} + \delta_{s+s'}\big)\big]^{\inv{2}},
\end{multline}
with $\bar p=p-2\geq 0$.

In parallel, a change $y\to x+y$ in \eqref{eq:scp_xy_first} provides 
\begin{equation*}
2\,\scp{J(x)}{y}\ \leq\ -\|x\|_p^2\ +\ (p-1)\,\|x+y\|_p^2 - \|y\|_p^2,
\end{equation*}
where we used the fact that $\scp{J(x)}{x}=\|x\|_p^2$. By summing this
inequality with (\ref{eq:scp_xy_first}), we have
\begin{equation*}
4\,\scp{J(x)}{y}\ \leq\ (p-2)\|y\|_p^2\ +\ (p-1)\,\|x+y\|_p^2 - \|x-y\|_p^2.
\end{equation*}
Using the RIP$_{p,2}$ on $x=\Phi u$ and $y=t\Phi v$ as above leads to
\begin{multline*}
4\mu_{p,2}^{-2}t\,|\scp{J(\Phi u)}{\Phi v}|\ \leq\ (1+\delta_{s'})\bar p\,t^2\\ 
+ (p-1)(1+\delta_{s+s'})(1+t^2) - (1-\delta_{s+s'})(1+t^2)\\
= \bar p + p \delta_{s+s'} + 
\big(2 \bar p + \bar p\delta_{s'}+p\delta_{s+s'}\big)t^2,
\end{multline*}
with the same argument as before to explain the absolute
value. Minimizing over $t$ as above gives 
\begin{multline} 
\label{eq-pl3-b}
2 \mu_{p,2}^{-2}\,|\scp{J(\Phi u)}{\Phi v}|\ \leq\\
\big[(\bar p + p\,\delta_{s+s'})\,\big(\,2\bar p +
\bar p\delta_{s'}+p\,\delta_{s+s'}\big)\big]^{\inv{2}}.
\end{multline}
Together, (\ref{eq-pl3-a}) and (\ref{eq-pl3-b}) imply
\begin{multline*}
C_p = \min \big\{\big[(\delta_{s}+\delta_{s+s'})\big(\delta_{s'} +
\delta_{s+s'} + \bar p\,(1 +
\delta_{s'}) \big)\big]^{\inv{2}},\\ 
\big[\big(\delta_{s+s'} + \bar p\,\tfrac{1 +
\delta_{s+s'}}{2}\big)\big(\delta_{s+s'} + \bar p\,\tfrac{2 +
\delta_{s'}+\delta_{s+s'}}{2}\big)\big]^{\inv{2}}\big\}.
\end{multline*}
It is easy to check that $C_p=C_p(\Phi,s,s')$ behaves as
$\sqrt{(\delta_{s}+\delta_{s+s'})\,(1+\delta_{s'})\,\bar p}$ for $\bar
p \gg \frac{\delta_{s'} + \delta_{s+s'}}{(1 + \delta_{s'})}$, and as
$\delta_{s+s'} + \tfrac{3}{4}(1+\delta_{s+s'})\bar p + O(\bar p^2)$
for $p\simeq 2$. 
\end{proof}

\section{}
\label{sec:proof-theor-refpr}

\begin{proof}[Proof of Theorem \ref{prop:l2-l1-instance-optimality-BPDQp}]
  Let us write $x^*_p=x+h$. We have to characterize the behavior of
  $\|h\|_2$. In the following, for any vector $u\in\Rbb^d$ with
  $d\in\{m,N\}$, we define $u_A$ as the vector in $\Rbb^d$ equal to
  $u$ on the index set $A \subset \{1,\,\cdots, d\}$ and 0 elsewhere.
  
  We define $T_0 = \supp x_K$ and a partition $\{T_k:1\leq k \leq
  \lceil (N-K)/K\,\rceil\}$ of the support of $h_{T_0^c}$. This
    partition is determined by ordering elements of $h$ off of the
    support of $x_K$ in decreasing absolute value. We have $|T_k|=K$
  for all $k\geq 1$, $T_k\,\cap\,T_{k'}=\emptyset$ for $k\neq k'$, and
  crucially that $|h_{j}|\leq |h_{i}|$ for all $j\in T_{k+1}$ and
  $i\in T_{k}$.

  We start from 
  \begin{equation}
    \label{eq-pt2-a}
    \|h\|_2\ \leq\ \|h_{T_{01}}\|_2\ +\ \|h_{T_{01}^c}\|_2,
  \end{equation}
  with $T_{01} = T_0 \cup T_1$, and we are going to bound separately the
  two terms of the RHS. In \cite{candes2008rip}, it is proved that
  \begin{equation}
    \label{eq:candes-compress-bound}
    \|h_{T_{01}^c}\|_2\leq \sum_{k\geq 2} \|h_{T_k}\|_2 \leq \|h_{T_{01}}\|_2\ +
    2e_0(K),
  \end{equation}
  with $e_o(K) = \tfrac{1}{\sqrt{K}} \norm{x_{T_0^c}}_1$.
Therefore, 
\begin{equation*}
  \|h\|_2\ \leq\ 2\|h_{T_{01}}\|_2\ + 2 e_0(K).
\end{equation*}

Let us bound now $\|h_{T_{01}}\|_2$ by using the
RIP$_{p,2}$. From the definition of the mapping $J$, we have 
\begin{multline*}
\|\Phi h_{T_{01}}\|_p^2\ =\ \scp{J(\Phi h_{T_{01}})}{\Phi
  h_{T_{01}}}\\
=\ \scp{J(\Phi h_{T_{01}})}{\Phi h}\ -\ \sum_{k\geq
  2}\,\scp{J(\Phi h_{T_{01}})}{\Phi h_{T_k}}.
\end{multline*}
By the H\"older inequality with $r=\frac{p}{p-1}$ and $s=p$,
\begin{multline*}
\scp{J(\Phi h_{T_{01}})}{\Phi h} \leq \|J(\Phi
h_{T_{01}})\|_r\|\Phi h\|_s\\ 
= \|\Phi h_{T_{01}}\|_p\|\Phi h\|_p\ \leq\ 2\,\epsilon\,\|\Phi
h_{T_{01}}\|_p\\
\leq  2\,\epsilon\,\mu_{p,2}\,(1+\delta_{2K})^{\frac{1}{2}}\|h_{T_{01}}\|_2,
\end{multline*}
since $\|\Phi h\|_p\leq \|\Phi x - y\|_p+ \|\Phi x^*_p - y\|_p\leq
2\epsilon$.  Using Lemma \ref{lemma:bound-scp-lp}, as $h_{T_{01}}$ is 2K
sparse and $h_{T_k}$ is K sparse, we know that, for $k\geq 2$,
\begin{equation*}
|\scp{J(\Phi h_{T_{01}})}{\Phi h_{T_k}}|\ \leq\
\mu_{p,2}^2\,C_p\,\|h_{T_{01}}\|_2\,\|h_{T_k}\|_2,
\end{equation*}
with $C_p=C_p(\Phi,2K,K)$, so that, using again the
RIP$_{p,2}$ of $\Phi$ and (\ref{eq:candes-compress-bound}),
\begin{multline*}
  (1-\delta_{2K})\mu_{p,2}^2\|h_{T_{01}}\|_2^2 \leq
   \|\Phi h_{T_{01}}\|_p^2\\[2mm]
  \leq
  2\epsilon\mu_{p,2}(1+\delta_{2K})^{\frac{1}{2}}\|h_{T_{01}}\|_2
   + \mu_{p,2}^2 C_p\|h_{T_{01}}\|_2\sum_{k\geq 2}\|h_{T_k}\|_2\\
\leq 2\epsilon\mu_{p,2}(1+\delta_{2K})^{\inv{2}}\|h_{T_{01}}\|_2\hfill\\
  +\mu_{p,2}^2C_p\|h_{T_{01}}\|_2\big(\|h_{T_{01}}\|_2 + 2e_0(K)\big).
\end{multline*}
After some simplifications, we get finally
\begin{equation*}
\|h\|_2\ \leq\
 \tfrac{2(C_p + 1-\delta_{2K})}{1-\delta_{2K} -
\ C_p}\ e_0(K) \ +\  \tfrac{4\sqrt{1+\delta_{2K}}}{1-\delta_{2K}
-\ 
 C_p}\ \tfrac{\epsilon}{\mu_{p,2}}.
\end{equation*}
\end{proof}

\section{}
\label{sec:proof-expec-and-val-lp-norm-unif-vec}
\begin{proof}[Proof of Lemma \ref{lemma:expec-and-val-lp-norm-unif-vec}]
  For a random variable $u\sim
  U([-\tfrac{\upalpha}{2},\tfrac{\upalpha}{2}])$, we compute easily that
  $\E|u|^p=\tfrac{\upalpha^p}{2^p(p+1)}$ and ${\rm
    Var}|u|^p=\tfrac{\upalpha^{2p}p^2}{2^{2p}(p+1)^2(2p+1)}$. Therefore,
  for a random vector $\xi\in\Rbb^m$ with components $\xi_i$
  independent and identically distributed as $u$, $\E\|\xi\|^p_p =
  \tfrac{\upalpha^p}{2^p(p+1)} m$ and ${\rm Var}\|\xi\|^p_p =
  \tfrac{\upalpha^{2p}p^2}{2^{2p}(p+1)^2(2p+1)}\, m$.

  To prove the probabilistic inequality below
  (\ref{eq:expec-and-val-lp-norm-unif-vec}), we define, for $1\leq
  i\leq m$, the positive random variables $Z_i =
  \tfrac{2^p}{\upalpha^p}|\xi_i|^p$ bounded on the interval $[0, 1]$
  with $\E Z_i = (p+1)^{-1}$. Denoting $S=\tinv{m}\sum_i Z_i$, the
  Chernoff-Hoeffding bound \cite{hoeffding1963pis} tells us that, for
  $t\geq 0$, $\Prob\big[S \geq (p+1)^{-1} + t\,\big] \leq e^{-2t^2 m}
  $.  Therefore,
\begin{equation*}
\Prob\big[ \|\xi\|_p^p\ \geq\ \tfrac{\upalpha^p}{2^p(p+1)}\,m + \tfrac{\upalpha^p}{2^p}\,t\,m \big] \ \leq\ e^{-2t^2 m},
\end{equation*}
which gives, for $t=\kappa\,m^{-\inv{2}}$,
\begin{equation*} 
\Prob\big[ \|\xi\|_p^p\ \geq\ \zeta_p\ +\ \tfrac{\upalpha^p}{2^p}\,\kappa\,m^{\inv{2}} \big] \ \leq\ e^{-2\kappa^2}.
\end{equation*}
The limit value of $(\zeta_p\ +\
\tfrac{\upalpha^p}{2^p}\,\kappa\,m^{\inv{2}})^{1/p}$ when $p\to\infty$
is left to the reader.
\end{proof}



\newpage
\begin{IEEEbiography}[{\includegraphics[width=25mm,keepaspectratio]{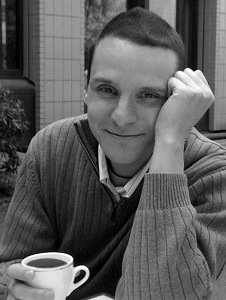}}]{Laurent
    Jacques} received the B.Sc. in Physics, the M.Sc. in Mathematical
  Physics and the PhD in Mathematical Physics from the Universit\'e
  catholique de Louvain (UCL), Belgium. He was a Postdoctoral
  Researcher with the Communications and Remote Sensing Laboratory of
  UCL in 2005-2006. He obtained in Oct. 2006 a four-year (3+1)
  Postdoctoral funding from the Belgian FRS-FNRS in the same lab. He
  was a visiting Postdoctoral Researcher, in spring 2007, at Rice
  University (DSP/ECE, Houston, TX, USA), and from 2007 to 2009, at
  the Swiss Federal Institute of Technology (LTS2/EPFL,
  Switzerland). His research focuses on Sparse Representations of
  signals (1-D, 2-D, sphere), Compressed Sensing, Inverse Problems,
  and Computer Vision.
\end{IEEEbiography}

\begin{IEEEbiography}[{\includegraphics[width=25mm,keepaspectratio]{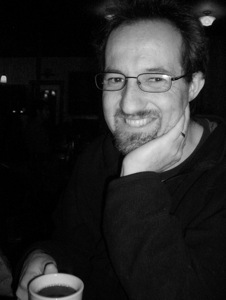}}]{David
    K. Hammond} was born in Loma Linda, California. He received a
  B.S. with honors in Mathematics and Chemistry from the Caltech in
  1999, then served as a Peace Corps volunteer teaching secondary
  mathematics in Malawi from 1999-2001. In 2001 he began studying at
  the Courant Institute of Mathematical Sciences at New York
  University, receiving a PhD in Mathematics in 2007. From 2007 to
  2009, he was a postdoctoral researcher at the Ecole Polytechnique
  Federale de Lausanne. Since 2009, he is postdoc at the
  NeuroInformatics Center at the University of Oregon, USA. His
  research interests focus on image processing and statistical signal
  models, data processing on graph, as well as inverse problems
  related to EEG source localization for neuroimaging.
\end{IEEEbiography}

\begin{IEEEbiography}[{\includegraphics[width=25mm,keepaspectratio]{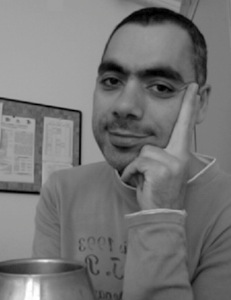}}]{Jalal
    M. Fadili} graduated from the Ecole Nationale Sup\'erieure
  d'Ing\'enieurs (ENSI) de Caen, Caen, France, and received the
  M.Sc. and Ph.D. degrees in signal and image processing from the
  University of Caen. He was a Research Associate with the University
  of Cambridge (MacDonnel-Pew Fellow), Cambridge, U.K., from 1999 to
  2000. He has been an Associate Professor of signal and image
  processing since September 2001 at ENSI. He was a visitor at several universities
  (QUT-Australia, Stanford University, CalTech, EPFL). His research interests
  include statistical approaches in signal and image processing,
  inverse problems, computational harmonic analysis, optimization and sparse
  representations. His areas of application include medical and astronomical imaging.
\end{IEEEbiography}

\end{document}